\documentclass[12pt]{iopart}

\usepackage{graphicx}
\usepackage{iopams}
\usepackage{amssymb}
\usepackage{amsthm}
\usepackage{stmaryrd}
\usepackage{bbm}
\usepackage{textcomp}

\usepackage{subfigure}

{\theoremstyle{plain}
	\newtheorem{thm}{Theorem}[section]
	
	\newtheorem{prop}[thm]{Proposition}

}

\newcommand{\SP}[0]{\mbox{\,$\mathcal{S}$}}

\newcommand{\Ct}{\mathcal{C}_2}
\newcommand{\Co}{\mathcal{C}_1}
\newcommand{\Gft}{F_2}
\newcommand{\Gfs}{F_{\rm S}}
\newcommand{\Gfo}{F_{\rm O}}
\newcommand{\Gfr}{F_{\rm R}}
\newcommand{\T}{\mathcal{I}}
\newcommand{\dd}[1]{\frac{\partial}{\partial #1}}
\newcommand{\EE}[3]{E\bigg[\frac{#1}{#2}{#3}\bigg]}
\newcommand{\lr}[1]{\left(#1\right)}

\begin{document}

\title{Exact generating function for 2-convex polygons.} 
\author{W. R. G. James, I. Jensen and A. J. Guttmann}

\address{
ARC Centre of Excellence for Mathematics and Statistics of Complex Systems, \\
Department of Mathematics and Statistics, 
The University of Melbourne, Victoria 3010, Australia}

\eads{\mailto{william.james@axa.com},\mailto{I.Jensen@ms.unimelb.edu.au},\mailto{T.Guttmann@ms.unimelb.edu.au}}

\date{\today}

\begin{abstract}
Polygons are described as almost-convex if their perimeter differs
from the perimeter of their minimum bounding rectangle by twice their
`concavity index', $m$. Such polygons are called \emph{$m$-convex}
polygons and are characterised by having up to $m$ indentations in
their perimeter. We first describe how we conjectured the (isotropic) generating function 
for the case $m=2$ using a numerical procedure based on series expansions.
We then proceed to prove this result for the more general case of
the full anisotropic generating function, in which steps in the $x$ and $y$ 
direction are distinguished. In so doing, we develop tools that would allow
for the case $m > 2$ to be studied. 
\end{abstract}

\submitto{\JPA}

\pacs{05.50.+q,05.70.Jk,02.10.Ox}

\maketitle

\section{Introduction} \label{s_intro}

\begin{figure}[h]
	\begin{center}
        \includegraphics[scale=1.0]{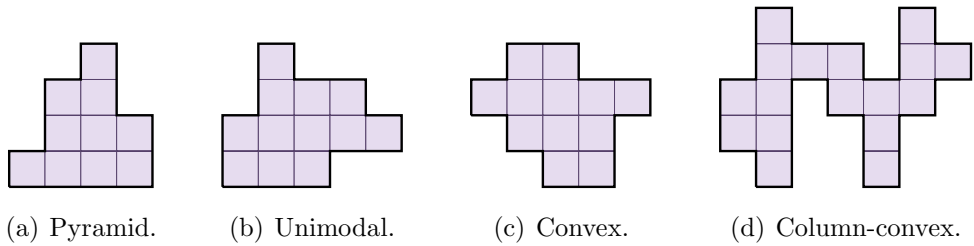}
	\end{center}
	\caption{Examples of polygons that are in some way convex.}
	\label{Fig:convex_examples}
\end{figure}

The enumeration of self-avoiding polygons (SAPs) is a classical problem in
statistical mechanics and combinatorics.
Exact results have thus far largely required the restriction of SAPs to sub-classes that
are in some way {\em convex}. In two dimensions, convexity means that the
perimeter is equal in length to the length of the minimum bounding rectangle (MBR).
Column-convexity means that any vertical cross-section may only intersect
the polygon twice, such that all columns are connected. Examples of convex and
column-convex polygons can be seen in figure~\ref{Fig:convex_examples}.
Convex polygons on two-dimensional lattices have been studied extensively by Lin
\cite{Lin:Chang88,Lin:90,Lin:91} and  Bousquet-M\'elou \cite{Bousquet:96} and
many exact results are known including the full area-perimeter generating function.
In 1997, Bousquet-M\'elou and Guttmann \cite{BMG:convex} gave exact
results for convex polygons in three dimensions and a method for their enumeration
in arbitrary dimensions \cite{BMG:convex}.

\begin{figure}
	\begin{center}
        \includegraphics[scale=1.0]{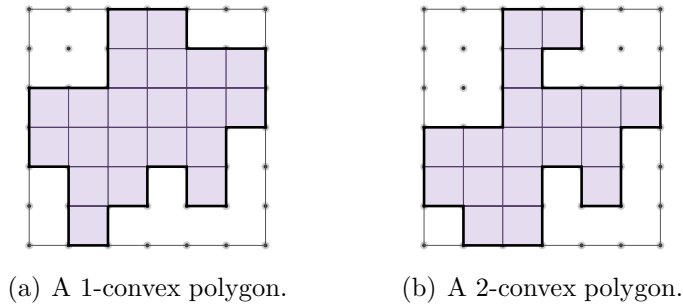}
	\end{center}
	\caption{Almost-convex polygons with the minimum bounding
	rectangle marked.}
	\label{Fig:almost_examples}
\end{figure}

Enting \etal \cite{Enting:92} described polygons as almost-convex if their
perimeter differs from the perimeter of their minimum bounding rectangle by
twice their `concavity index', $m$. Such polygons are called \emph{$m$-convex}
polygons and are characterised by having up to $m$ indentations in their perimeter.
Examples of 1-convex and 2-convex polygons can be found in
figure~\ref{Fig:almost_examples}.
Enting \etal derived the asymptotic behaviour of the
number of $m$-convex polygons according to their perimeter, $n$ for $m =
o(\sqrt{n})$. The results were confirmed for the case $m=0$ (i.e. convex
polygons) by the known perimeter generating function.  Subsequently,
 Lin~\cite{Lin:1-convex} derived the exact generating function for
1-convex polygons, using a `divide and conquer' technique introduced to the
problem of convex animals (the interior of a convex SAP) by Klarner and Rivest
\cite{Klarner:74}. His result provided support for a conjecture in
\cite{Enting:92}, giving the next term in the asymptotic expansion for the number
of polygons with perimeter $n$ and concavity index $m.$

This is the second in a series of papers that look at families of $m$-convex
polygons. In the first \cite{wrgj} we outlined the 50-year history of polygon 
enumeration on the square lattice before re-deriving the generating functions for
1-convex polygons in an effort to generalise the methodology and extend the
results to osculating\footnote[1]{Osculating polygons are those that may touch
themselves, but not cross.} and neighbour-avoiding\footnote[2]{Neighbour-avoiding
polygons are those that may not occupy a neighbouring lattice vertex without
being connected by an edge.} polygons.

Polygon models have long been used to model vesicles, with
self-avoiding polygons being the canonical model \cite{Fisher:91}. Associating a
fugacity with the {\it area} of the polygon, a phase transition
occurs, so that at sufficiently high fugacity, the polygons change
from the self-avoiding class to the convex class. The introduction of
$m$-convex polygons permits the exploration of this phase transition
in greater detail. Firstly, we find the (unsurprising) result that
$m$-convex polygons, for $m$ finite, have the same fractal dimension 
(and hence are in the same universality class) as convex polygons,
that is to say, the fractal dimension remains unchanged at 2
(as compared to the value $4/3$ for self-avoiding polygons).
Secondly, the $m$-convex model permits one to
associate a fugacity with the concavity index $m$, and this would
be the polygon analogue of the stiffness in self-avoiding walk models
of polymer stretching. It is however not our purpose to study this aspect of
the problem here. Ideally, one would like to predict how the form of the 
generating functions changes as the concavity index grows. This would help us 
understand what happens in the scaling limit as the concavity index grows 
in proportion to the perimeter.

In 2000, in unpublished work, we conjectured, on the basis of long series expansions, the
(isotropic) generating function for the case $m=2.$ In this paper, we describe 
the numerical procedure that led to this conjecture, and proceed to prove it. Indeed, 
we do so for the full anisotropic generating function, in which steps in the $x$ and $y$ 
directions are distinguished. In so doing, we develop tools that allow
for the case $m > 2$ to be studied, though we do not do so.

In section~\ref{s_enum} we describe the series expansions that allowed us to conjecture
the exact result, and also, not incidentally, to provide checks on our rigorous results in
the process of proving our conjecture.
In the following section we introduce the
methodology used to derive our results, followed by examples of its application.
We enumerate all the separate
building blocks required in the factorisation of 2-convex polygons.
The intermediate results, as well as much of the detail, are omitted for reasons of
conciseness and clarity. These may be found together with the equivalent staircase and
unimodal results in \cite{thesis} and \cite{JJG07}, where the presented results first
appeared.

\section{Definitions and notation} \label{s_defs}

In this section we briefly summarise some definitions and notation used in the remainder
of the paper.

\begin{figure}
\begin{center}
\includegraphics[scale=0.5]{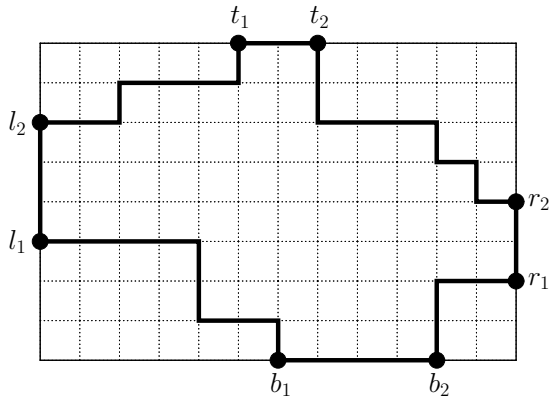}
\end{center}
\caption{A convex polygon with marked vertices where the polygon meets
the minimum bounding rectangle }
\label{Fig:Convex_def}
\end{figure}

\begin{description}
\item{Directed Walks (DWs):} Square lattice walks that take either positive or
negative steps in each of the horizontal and vertical directions, but not both.
For example, a walk that only steps up and to the left is directed.
\item{Generating functions:} If the number of polygons (in a given class) 
with perimeter $n$ is $p_n$ then the associated isotropic generating function 
is $F(x)=\sum_n p_nx^n$.  For polygons on the square lattice the perimeter is 
always even and we shall therefore study the {\em half-perimeter} generating
functions  $G(x)=\sum_n p_{2n}x^n$.
In more generality we distinguish between the number of steps in the $x$ and $y$
directions and study the full anisotropic generation functions,
$G(x,y)=\sum_{m,n} p_{2m,2n}x^m y^n$, where $p_{2m,2n}$ is the number of polygons with
$2m$ horizontal steps and $2n$ vertical steps. 
\item{The half-perimeter operator:} We denote by $E$ the operator that converts the
perimeter generating function to the half-perimeter generating function. $E_x$
(resp. $E_y$) converts only the direction counted by $x$ (resp. $y$). (For a
full definition, see \cite{BMG:convex}.)
We therefore have
\begin{equation}
	E\left[f(x,y)\right] = E_x\left[E_y\left[f(x,y)\right]\right]
\end{equation}
where
\begin{equation}
	E_x\left[f(x,y)\right] = (f(\sqrt{x},y)+f(-\sqrt{x},y))/2.
\end{equation}
If $x$ or $y$ (or a function of $x$ or $y$) is asterisked, then the operator only
takes the half-perimeter of the non-asterisked generating function. That is,
\begin{equation}
	E_x\left[f(x,x^*)\right] = E_x\left[f(x,y)\right]\big|_{y\rightarrow x}.
\end{equation}
For example, $E_x[1/(1-x-x^*)] = (1-x)/(1-3x+x^2)$ and, more generally,
\begin{equation}
	E_x[f(x)g(x^*)] = g(x) E_x[f(x)],
\end{equation}
which extends naturally to the multivariate case.
This allows the conversion from the perimeter generating function to the half-perimeter
generating function for different factors of a polygon separately.
\item{Minimum bounding rectangle (MBR):} This is the smallest rectangle which encloses the
polygon. The polygons in figure~\ref{Fig:almost_examples} have a $6\times6$ MBR while the
polygon in figure~\ref{Fig:Convex_def} has a $12\times 8$ MBR.
\item{Sides:} A convex polygon can be factored, as per figure~\ref{Fig:Convex_def},
into four overlapping DWs: from $l_1$ to $t_2$, from $t_1$ to $r_1$,
and so on. These are the maximal directed factors of the polygon, and we refer
to them as the {\it sides} of the polygon.
\item{Arcs:} An {\it arc} of a convex polygon is made up of a pair of adjacent sides.
That is, it a maximal partially-directed factor of the polygon.
Referring to figure~\ref{Fig:Convex_def}, the top {\em arc} is the
path from $l_1$ to $r_1$, passing through  $ l_2, t_1, t_2 $ and $ r_2 $,
and similarly for the bottom, left, and right arcs.
\item{Indents:} An indent occurs in a side when the DW takes a step
in the `wrong' direction. We refer to such walks as `almost-directed'.
As an example, consider the part of the perimeter
 on the top-right side, which, when the polygon is traversed anti-clockwise, 
only takes steps up and to the left. An indent would occur if this walk were to take
steps down (or to the right) and then resume taking steps up and to the left. 
We shall distinguish between indents in the vertical and horizontal directions.
In figure~\ref{Fig:almost_examples}a the polygon has a single vertical indent while the polygon
in  figure.~\ref{Fig:almost_examples}b has both a vertical and horizontal indent. 
A walk taking $k'$ steps in the `wrong' direction followed
by $k''$ steps in the `correct' direction produces an indent $k=\min(k',k'')$ deep.
Note that a vertical indent may contain a horizontal indent and vice versa.
\end{description}

\section{Exact solution from enumerations} \label{s_enum}

Several years ago two of us (IJ and AJG) found the exact generating
function, $\Ct(x)$ for 2-convex polygons numerically from exact enumerations
for 2-convex polygons and some simpler sub-classes. We found the solution
by counting the number of 2-convex polygons by using a program
designed to enumerate SAPs \cite{Jensen99}. This program counts
the number of SAPs by perimeter having a given MBR. From this data
it is trivial to extract the 2-convex data. However, with the computational
resources available at the time, we could not get a series long enough
to find the generating function directly (we counted 2-convex polygons up to
perimeter 110 yielding 48 non-zero terms). We therefore enumerated
three sub-classes, namely 2-convex polygons with one vertical indent 2 deep
on the top arc, 2-convex polygons with  two vertical indents 1 deep
on the top arc, and 2-convex polygons with vertical indents 1 deep
on the top and bottom arcs. The programs used in these enumerations were
simple generalisations of the one used by Guttmann and Enting \cite{Guttmann:Enting88}
to count convex polygons. We then used the series for these sub-classes
to find the respective generating functions $\Gft(x)$, $\Gfs(x)$ and $\Gfo(x)$.
Finally, we looked at the `remainder' of the full 2-convex case
$\Gfr(x) = \Ct(x)-4\Gft(x)-4\Gfs(x)-2\Gfo(x)$, which counts
cases of 2-convex polygons with vertical and horizontal indents 1 deep,
and managed to find the generating function (and hence the full 2-convex
generating function). 

We found the solution to the special cases via a judicious guess for the
form of the generating functions. From Lin \cite{Lin:1-convex}, the isotropic
generating function for 1-convex polygons is:

\begin{eqnarray}\label{eq:1Cgf}
\fl
\qquad \Co(x) =  {\frac {{x}^{3} \left( -4+56x-300{x}^{2}+773{x}^{3}-973{x}^{4}
+535{x}^{5}-90{x}^{6}+24{x}^{7} \right) }{ \left( 1-x \right) 
 \left( 1-3x+{x}^{2} \right)  \left( 1-4x \right) ^{3}}} \nonumber \\
+{\frac {4{x}^{3} \left( 1-9x+25{x}^{2}-23{x}^{3}+3{x}^{4}
 \right) }{ \left( 1-x \right)  \left( 1-4x \right) ^{5/2}}}. 
\end{eqnarray}
From this it is reasonable to expect that the generating function for 2-convex polygons 
$\Ct(x)$ and the special cases $\Gft(x)$ and so on are of a similar form,
$[A(x)+B(x)\sqrt{1-4x}]/D(x)$, where $A(x)$, $B(x)$ and $D(x)$ are polynomials.
In particular we expect the denominator $D(x)$ to be similar to the one
in the above expression, but with larger exponents and possibly involving
further simple factors. If we can find $D(x)$ then $A(x)$ and $B(x)$ can
be found simply from a formal series expansion using polynomials
with unknown coefficients. By equating the terms in this formal series 
with the known series for say $\Gft(x)$ we find a set of linear equations
for the unknown polynomial coefficients. 

Here we give some further details of how we found the generating function $\Gft(x)$.
We calculated the number of convex polygons with an indent 2 deep on the top arc
to perimeter 216. This gives us the first 100 non-zero terms in the half-perimeter
generating function. Our first task is to determine the denominator $D(x)$. We
did this by analysing the singularities of  $\Gft(x)$ using differential 
approximants. Our analysis showed that the series has singularities at 
$x=1/4$ with exponents $-3$ and $-5/2$ (this confirms that at the 
dominant singularity we have a square-root correction term), at $x=0.381966\ldots$ (the
first root of $1-3x+x^2$) with exponents $-3$ and $-1$, at $x=1$ with
exponents $-4.99(3)$ and $-3.0(5)$, and at $x=2.618\ldots$ 
(the second root of $1-3x+x^2$) with exponent $-3$. The conclusion 
is that in this case $D(x)=(1-x)^5(1-3x+x^2)^3(1-4x)^3$. 
By inserting this into the general form and equating terms in the
formal expansion with those of $\Gft(x)$ we found a solution with polynomials
$A(x)$ and $B(x)$ of degree 18 and 17, respectively, with $B(x)$ containing
the factors $(1-x)^2(1-3x+x^2)^2$, the latter of these factors was
indicated by the exponents found at $x=0.381966\ldots$. The polynomials
are:
\begin{eqnarray}\label{eq:2Deeppoly}
\fl
\qquad A(x)=-8{x}^{2}+208{x}^{3}-2428{x}^{4}+16856{x}^{5}-77742{x}^{6}+252114{x}^{7} \nonumber \\
-593563{x}^{8}+1032521{x}^{9}-1336471{x}^{10}+1284072{x}^{11}-904540{x}^{12}  \nonumber \\
+456064{x}^{13}-158327{x}^{14}+36093{x}^{15}-4955{x}^{16}+126{x}^{17}+88{x}^{18}. \nonumber \\
\fl
\qquad B(x)=(1-x)^2(1-3x+x^2)^2(8{x}^{2}-128{x}^{3}+844{x}^{4}-2992{x}^{5} \nonumber \\
+6262{x}^{6}-8014{x}^{7}+6188{x}^{8}-2602{x}^{9}+470{x}^{10}-12{x}^{11}).
\end{eqnarray}

Similarly we found the generating functions for the other three special cases
$\Gfs(x)$, $\Gfo(x)$ and $\Gfr(x)$. The only additional point worth noting is that
the denominator in the case of $\Gfr(x)$ contains the extra factor $(1-2x)$.
Collating these results we find that:

\begin{equation}\label{eq:2Cgf}
\fl
\qquad \Ct(x)= \frac{A_2(x)}{(1-x)^7(1-2x)(1-3x+x^2)^3(1-4x)^4}+\frac{B_2(x)}{(1-x)^3(1-4x)^{7/2}},
\end{equation}
where
\begin{eqnarray}\label{eq:2Cpoly}
\fl
\qquad A_2(x)= -24{x}^{2}+864{x}^{3}-14368{x}^{4}+146672{x}^{5}-1030216{x}^{6}+5289512{x}^{7} \nonumber \\
-20587766{x}^{8}+62176564{x}^{9}-147946110{x}^{10}+280112802{x}^{11} \nonumber \\ 
-424512212{x}^{12}+516373058{x}^{13}-504068274{x}^{14}+393649476{x}^{15} \nonumber \\
-244279626{x}^{16}+119050550{x}^{17}-44773540{x}^{18}+12722814{x}^{19} \nonumber \\
-2660520{x}^{20}+378184{x}^{21}-22560{x}^{22}-3200{x}^{23}+512{x}^{24}.\nonumber \\
\fl
\qquad B_2(x)=-24{x}^{2}+456{x}^{3}-3592{x}^{4}+15264{x}^{5}-38200{x}^{6}+57792{x}^{7} \nonumber \\
-52832{x}^{8}+28872{x}^{9}-8968{x}^{10}+1248{x}^{11}+128{x}^{12}.
\end{eqnarray}

In the next sections we show how to prove this result for the general
anisotropic case. 

\section{Enumeration techniques}

\subsection{Convex polygon basics}
Following from Section~\ref{s_defs}, one may describe convex polygons as a series
of four non-intersecting DWs that make up the four sides of the polygon. (To see
this, we refer to figure~\ref{Fig:Convex_def}.) If one of these sides has no
steps in the interior of the MBR, this means that the polygon touches one of the
corners of the MBR. This sub-class of convex polygons is referred to as
directed-convex\footnote{
	This name comes from the definition which says that all cells in the
	interior of the polygon can be connected with the corner cell by a directed
	walk on the dual graph.
} or unimodal\footnote{
	This name comes from the fact that there is only one mode in each direction,
	when we take the projection of the walk in that direction.
}. We could therefore define this class of polygons as having only three sides.
Similarly, staircase polygons are formed of two DWs that start and end at diagonally
opposite corners of the MBR. Finally, pyramids and stack polygons also have two sides
formed by DWs, but these are adjacent sides, with a straight base or side edge.
(Stack polygons are simply pyramids on their sides.) Such classes of convex
polygons were depicted in figure~\ref{Fig:convex_examples}.

For notational convenience, let us now define some well-known generating
functions, where $x$ (resp. $y$) counts the horizontal (resp. vertical) steps.
The generating function for pairs of intersecting DWs that
begin and end at the same points (referred to as staircase {\em festoons}) we
denote as
\begin{equation}
	Z = \sum_{n,m} {{n+m}\choose{n}}^2 x^n y^m = 1/\sqrt{\Delta},
\end{equation}
where $\Delta=1-2x-2y-2xy+x^2+y^2$. The staircase polygon generating
function we denote as
\begin{equation}
	\SP = (1-x-y-\sqrt{\Delta})/2.
\end{equation}
We note that the unimodal generating function is simply $xyZ$. Now, by defining
\begin{equation}
	u = x + \SP \quad \mbox{and} \quad v = y + \SP,
\end{equation}
we can re-express all our almost-convex polygon generating functions as
expressions with terms that are simply the quotient of polynomial functions
of $u$ and $v$. This is achieved via the transformation of variables
\begin{equation}
	x = u(1-v) \quad \mbox{and} \quad y = v(1-u).
\end{equation}
For example, we have
\begin{equation}
	\Delta = (1-u-v)^2,\quad  Z = 1/(1-u-v)\quad \mbox{and} \quad \SP = uv.
\end{equation}

\subsection{Joining polygon factors}
The Temperley method is central to the enumeration of partially convex polygons.
The so-called `functional-Temperley' method allowed Bousquet-M\'elou \cite{Bousquet:96} 
to enumerate classes of column-convex polygons.  It differs from the Temperley method 
in that it can be used to concatenate several large enumerable parts of the polygon, 
rather than individual columns. In particular the concatenated building blocks may 
be different types of polygons.  A variation \cite{BoRe:02}
allowed for the enumeration of certain classes of animals, represented as heaps
of dimers. Rechnitzer \cite{Rechnitzer} identified these methods as equivalent,
the superiority of one over the other lying in its ease of use and
appropriateness to the recurrence relation underlying the problem.

We call distinct enumerable parts of the polygons {\em factors}, due to our
ability to factorise the polygons into such parts by separating them at
unique factorisation points. In this paper, we use {\em factorisation
lines}, which we define by extending the interior edge of the indentations of
almost-convex polygons into lines that bisect the lattice (see
figure~\ref{Fig:1-convex}.)

\subsubsection{The Hadamard product} \label{s_Hadamard}

The approach that has been used repeatedly by Lin
\cite{Lin:Chang88,Lin:90,Lin:91} in the enumeration of convex polygons is to
build them up vertically, block by block.  The functional-Temperley method can
therefore be used in this case.  One tool which can be used to `join' polygonal
blocks together is the Hadamard product. (For a full description, see
\cite{Rechnitzer}.)  This is particularly useful when a few blocks need to be
joined in a non-recurring manner. We use the Hadamard product in calculating
most of the generating functions in this paper.

	Consider two series, $f(t)=\sum_n f_n t^n$ and $g(t)=\sum_n g_n t^n$,
	denote by $\odot_t$ the \emph{Hadamard product} with respect to $t$, then
	\begin{equation} \label{Eq:hadamard}
		f(t) \odot_t g(t) = \sum_n f_n g_n t^n.
	\end{equation}
	The \emph{restricted Hadamard product} with respect to $t$ is defined as
	\begin{equation} \label{Eq:r-hadamard}
		f(t) \odot_t g(t) \arrowvert_{t=1} = \sum_n f_n g_n
		= \frac{1}{2\pi i} \oint f(t) g(1/t) \frac{\mbox{d}t}{t}.
	\end{equation}
	For notational convenience, we will refer to this as a \emph{Hadamard
	join} (over $t$), or simply a `join'.

\begin{figure}
	\begin{center}
		\includegraphics{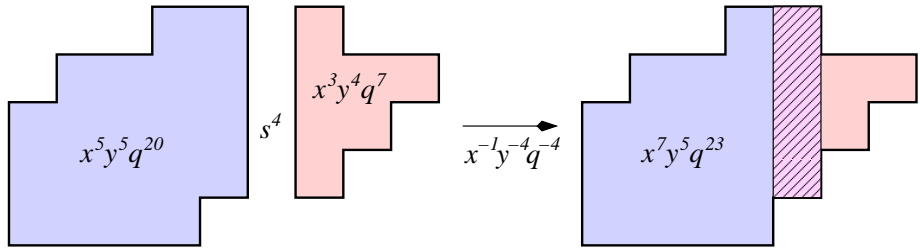}
	\end{center}
	\caption{An example of the action of the Hadamard product.}
	\label{Fig:hadamard}
\end{figure}

Hence, the Hadamard product over a given variable, say $s$, is the operator
which `joins' generating functions by matching the perimeters enumerated by
$s$. This is equivalent to joining the polygons by matching the edges of the
respective polygons, such that they overlap, and then removing the overlapping
edges, forming a single, larger polygon. We can also match the edge column,
rather than the edge perimeter. For example, if we enumerate
staircase and stack polygons according to their right perimeter, total
perimeter and area, we can join them, making the neighbouring columns overlap.
Making the transformation $s\mapsto s/yq$ and dividing by $x$, so that the
overlapping column is not double-weighted, we form unimodal polygons, as in
figure~\ref{Fig:hadamard}.

\subsubsection{Hadamard arithmetic}

It is straight-forward to show the following properties of restricted Hadamard
products \cite{Rechnitzer}.

\paragraph{It is distributive:}
\begin{equation}
	f(t) \odot_t (g(t)+h(t)) = f(t) \odot_t g(t) + f(t) \odot_t h(t).
\end{equation}

\paragraph{It follows the product rule:}
\begin{equation}
	\dd s (f(s,t) \odot_t g(s,t)) = ( \dd s f(s,t)) \odot_t g(s,t).
	+ f(s,t) \odot_t ( \dd s g(s,t) )
\end{equation}

\paragraph{It evaluates simply at poles:}
\begin{eqnarray}
	f(t) \odot_t \frac{1}{1-\alpha t} & = f(\alpha),
\\	f(t) \odot_t \frac{t^k k!}{(1-\alpha t)^{k+1}} & = \left(\dd t\right)^k f(t)
	\bigg|_{t=\alpha}.
\end{eqnarray}

We showed in \cite{wrgj} how $u$ and $v$ can be used to simply express the
generating functions for staircase polygons with fixed steps in the corner.
The generating function $u^av^b$ counts staircase polygons that start with $a$
horizontal steps and end with $b$ vertical ones (or {\em vice versa}). This
allows us to enumerate parts of polygons with specific sequences of steps along
the joins. Moreover, because we can express the generating functions for all
factors of almost-convex polygons as the quotient of polynomial expressions of
$u$ and $v$, evaluating the joins generally becomes straight-forward by
separating the poles in the denominator using partial fractions.

\subsection{Distinguishing steps to insert indents}
\begin{figure}
	\begin{center} \includegraphics[scale=0.8]{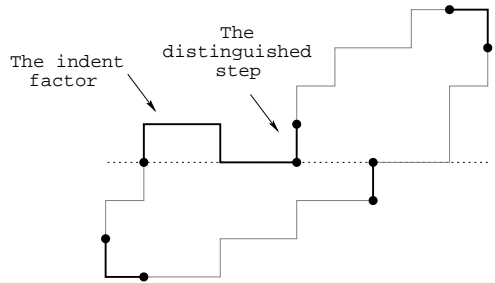} \end{center}
	\caption{The form of a 1-staircase polygon, which may be formed by inserting
	an indent factor next to a distinguished step.}
	\label{Fig:1-stair}
\end{figure}

One way of inserting indents in convex polygons is to distinguish a step for the
location of the indent and then make the appropriate adjustment to the
generating function. We therefore factorise the polygon at the distinguished
step by extending a line perpendicular to it. 
When an indent is joined to a staircase factor (as per figure~\ref{Fig:1-stair}),
the adjustment required for the indent is independent of its location.

We distinguish between the walk forming the indent and the rest of
that side of the polygon. In the case of a vertical indent
(as depicted in the figure), the indent starts with the vertical step at 
the same height as the distinguished step preceding the non-directed vertical step(s), 
and includes all steps up to (but not including) the distinguished step.
We refer to it as the {\em indent factor}.  We note that a single, $m$-deep
indent factor is therefore in the form of a pyramid.

As mentioned above, the generating function for staircase polygons with $a$
(resp. $b$) fixed horizontal (resp. vertical) steps in one corner is $u^av^b$.
As each extra fixed step along the factorisation line contributes $u$
to the generating function, the contribution of the indent factor to these
otherwise staircase polygons (which are called 1-staircase polygons) is
therefore $u^2/(1-u)^2$.
We therefore define the 1-deep indent generating function to be
\begin{equation} \label{Eq:I}
	\T = \frac{u^2}{(1-u)^2}, \quad \mbox{with} \quad
	\bar{\T} = \T(y,x) = \frac{v^2}{(1-v)^2},
\end{equation}
and, for $m>1$, the $m$-deep indent generating function is
\begin{equation} \label{Eq:Im}
\T_m \equiv \T_m(x,y) = \frac{u^2}{(1-u)^{2m}} = \frac{v^2 \SP^{2m-2}}{y^{2m}},
	\quad \mbox{with} \quad \bar{\T}_m = \T_m(y,x).
\end{equation}

The height-independent insertion of indents can be extended to unimodal polygons.
Whereas staircase
polygons can be separated into two halves -- one with only positive steps, the
other with only negative steps -- unimodal polygons are defined by the fact that
their positive horizontal steps occur before the negative ones, and similarly
for the vertical steps. This means that if the unimodal polygon
intersects itself, a staircase factor is formed. And so, if we factor a
unimodal polygon along the base of a vertical indent on the left arc (as per the
staircase factorisation shown in figure~\ref{Fig:1-stair}), the indent
must be a part of a staircase factor to the bottom-left and is therefore
enumerated by $\T_m$. This leads us to the following proposition.
\begin{prop}
	The generating function for bimodal $m$-staircase (resp. $m$-unimodal)
	polygons that are rooted in the bottom left corner and whose single
	$m$-deep indent is vertical and on the left arc is
	\(
		\T_m y^2 \dd y  (\mathcal{Q}/y),
	\)
	where $\mathcal{Q} = \SP/y$ (resp. $xZ$).
\end{prop}

This argument can be extended to enumerate unimodal polygons with two distinct
indents on the same side.  To construct such polygons, one may insert both
indents at a distinguished height, and then mark a second height where we
would like the second indent. We can therefore try to form the desired polygons by
translating the closest of the two indents to the second of the distinguished
heights. In \cite{JJG07}, we
point out that if the second indent lies below the first, this downward
translation may cause the polygon to intersect. Furthermore, in the other case,
an upward translation will mean that the bottom arc of the polygon will not go
above the original position of the translated indent. We show, however, that the
missing polygons in one case are equal in number to the extra polygons in the
other case. This intriguing fact leads to the following proposition.
\begin{figure}
	\begin{center}
	\subfigure[The form of 2-staircase polygons with both indents on their top-left side.]
		{\qquad\includegraphics[scale=0.4]{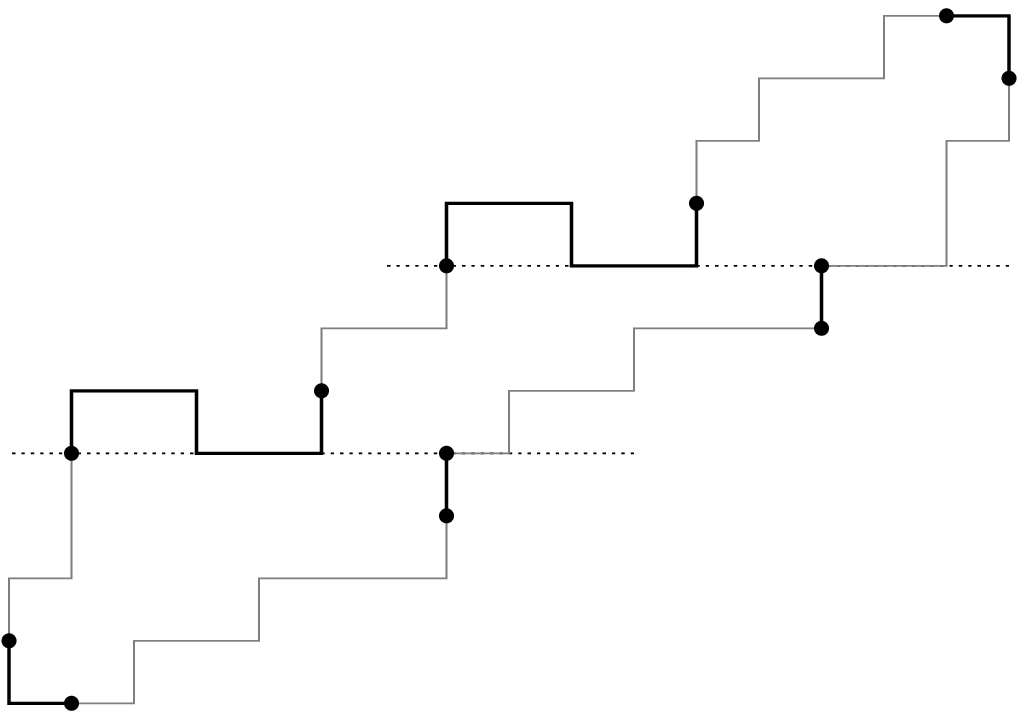} \qquad} \quad
	\subfigure[The two possible ways of inserting two indents.]
		{\includegraphics[scale=0.4]{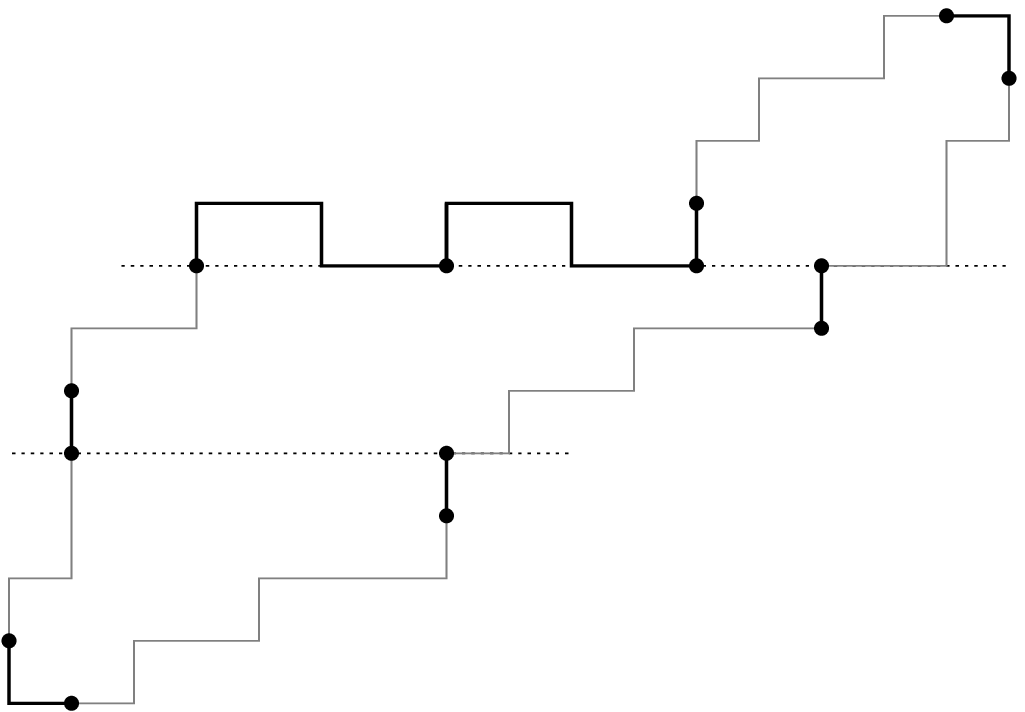} 
		 \includegraphics[scale=0.4]{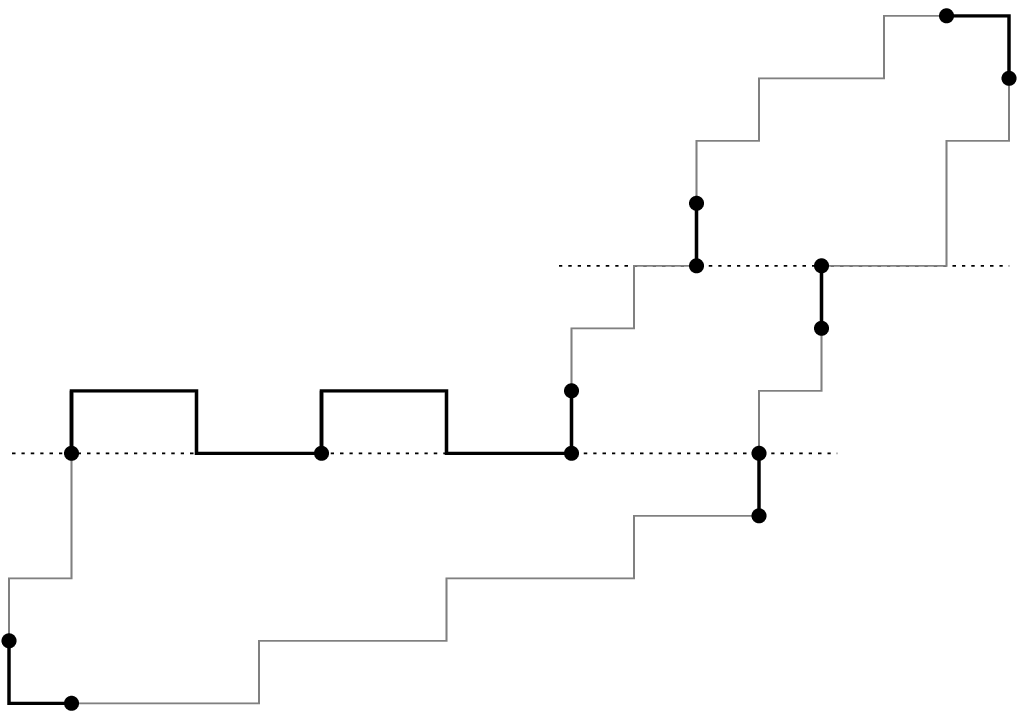} }
	\end{center}
	\caption{The form of 2-staircase polygons with both indents on their top-left side.}
	\label{Fig:2-stair}
\end{figure}

\begin{prop}
	The generating function for 2-staircase (resp. 2-unimodal) polygons that are
	rooted in the bottom left corner and with
	distinct vertical indents on the left arc is
	\begin{equation}
		y^3 \dd y (\T^2 \dd y (\mathcal{Q}/y))/2,
	\mbox{\quad where\quad} \mathcal{Q} = \SP/y \ (\mbox{resp. } xZ).
	\end{equation}
\end{prop}

\subsection{Folding Walks and Wrapping Polygons}

\subsubsection{Folding Walks}
An important notion in deriving exact generating functions for almost-convex
polygons is that of the so-called `folding' of DWs. This notion is simply a way
of describing the combinatorial objects enumerated by the enumerative methods
developed in \cite{BMG:convex}. These methods use the half-perimeter operator
defined in Section~\ref{s_defs} to enumerate intersecting convex polygons, and
then remove those that intersect. We are therefore interested in the enumeration
of these intersecting polygons.

We begin by considering DWs, which are enumerated by
\(
	1/(1-x-y),
\)
where $x$ (resp. $y$) counts the horizontal (resp. vertical steps). If we then
only consider those with an even number of horizontal steps, we can find either
a vertex or a series of vertical steps that have half of these horizontal steps
to the left, and half to the right. We can then fold the walk vertically at this
point by reflecting all the steps to the left over to the right. (The axis of
reflection is the vertical line that goes through the half-way point.) This walk
is now half as wide as it was, and its generating function is therefore
\begin{equation}
	E_x\left[\frac 1{1-x-y}\right].
\end{equation}
By removing those that have a horizontal step after the fold, we obtain
self-avoiding walks, which is an example of standard inclusion-exclusion
techniques,
\begin{equation}
	E_x\left[\frac {1-x}{1-x-y}\right].
\end{equation}
Adding a width-one column to its side then gives us stack polygons (that is,
sideways pyramids).

Now, folding vertically as well allows us to make the walk end at its
origin, forming a polygon. Forcing the polygon to start with a horizontal step and end
with a vertical one (or {\em vice versa}) then gives the following generating
function for (possibly intersecting) unimodal polygons:
\begin{equation}
	E\left[\frac {xy(1-x)(1-y)}{1-x-y}\right].
\end{equation}

Finally, by distinguishing a horizontal step after which we fold (rather than
choosing the half-way point) forces a second fold on the other side of the
walk in order to make its width equal to half its original width. This means
that the resulting polygon does not necessarily visit any corner of the MBR, and
we are left with (possibly intersecting) convex polygons. The resulting
generating function is
\begin{equation}
	E\left[\frac {xy(1-x)^2(1-y)^2}{(1-x-y)^2}\right],
\end{equation}
which is the $d=2$ case of Bousquet-M\'elou and Guttmann's multi-dimensional
result (Lemma 2.2, \cite{BMG:convex}).

Folding walks is therefore a simple way of enumerating intersecting polygons (or
factors of polygons, for that matter) with convexity requirements. This will
allow us, in the following section, to enumerate 1-convex polygons in a direct,
closed-form expression that is combinatorially interpretable. This provides us
with an example of the methods required for the enumeration of the various
sub-classes of 2-convex polygons. However, first we need to be able to fold one
factor within a polygon (that is joined directly to another factor) without our
methods breaking down. And for this we need ``wrapping''...

\subsubsection{Wrapping Polygons}

`Wrapping' refers to folding a single factor of a polygon that has been
constructed by joining multiple factors together.
As an example consider 1-unimodal polygons formed by joining a
staircase bottom factor to a unimodal top factor and an indent factor. The
total height of the polygon is given
by the sum of the heights of the top and bottom factors. The total width of the
polygon is measured by adding the width of the top factor to the width of the
bottom factor \emph{that lies to the left of the top factor}. In the 1-unimodal
case, the generating function for the bottom factor is therefore $\T v(u/x)^d$,
where $d$ is the length of the join.

\begin{figure}
	\begin{center}
	\subfigure[The form of the folded top factor.]
				{\qquad\includegraphics[scale=0.8]{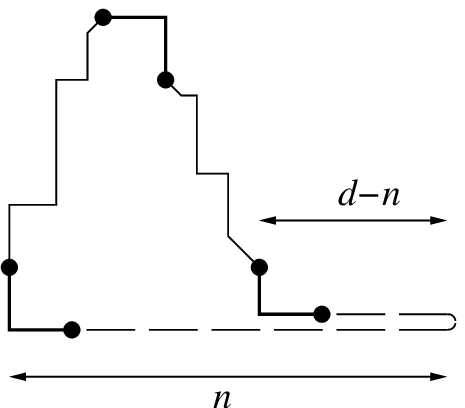}\qquad} \qquad
	\subfigure[The form of the wrapped polygon.]
			{\raisebox{10pt}
			{\hspace{30pt}
				 \includegraphics[scale=0.4]{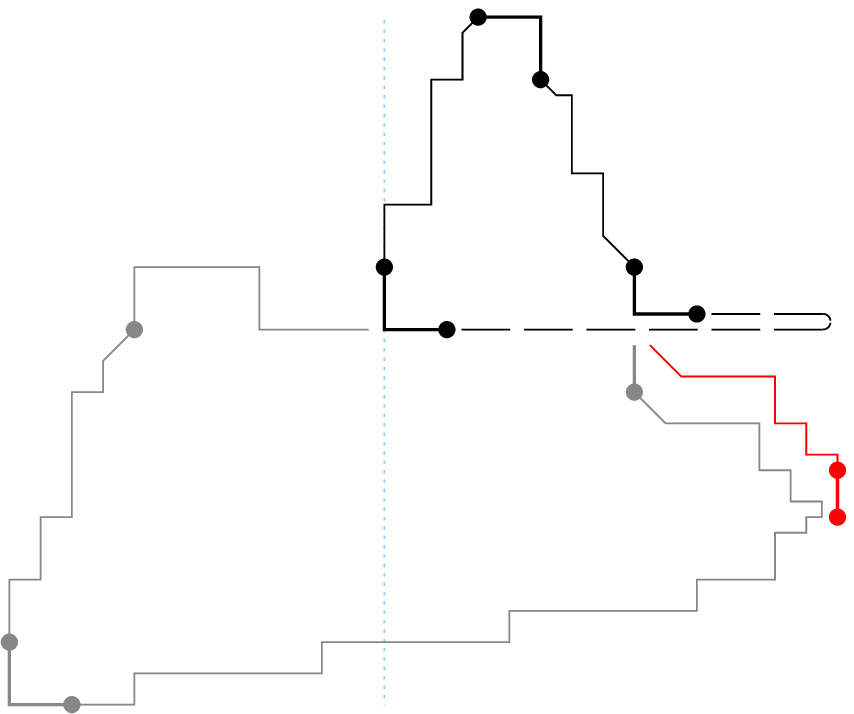}\hspace{40pt}}}
	\end{center}
	\caption{The action of wrapping when the top factor is narrower than the
	length of the join, $d$.}
	\label{Fig:wrapping}
\end{figure}

In this case, the unimodal factor is enumerated as a folded walk, as per the
previous section. However, we will sometimes fold the fixed steps of the
polygon, leaving  a chain of double-bonds of fixed steps around the fold, as
shown in figure~\ref{Fig:wrapping}(a).
If the join is of length $d$, then there are $d$ fixed horizontal steps in
each of the top and bottom factors that are identified, but are then removed and
do not form part of the polygon.
And so, if there are fewer than $d$ horizontal steps in the rest of the top
unimodal factor, as depicted in the figure, then the contribution
to the polygon is a pyramid of width $2n-d$, with a weight of $x^n$.
Importantly, although those fixed steps are not part of the polygon itself,
they do contribute to the weight. This is because the width of the polygon here
is given by the top factor. When the $E$ operator folds the top factor, in order
to keep the bottom factor joined to the top factor, it must therefore also be
folded. We say that we have {\em wrapped} the bottom factor. We can therefore see
that the fixed steps along the join that were folded to form double-bonds and
whose width is counted in the resulting generating function is the
projection of the wrapped part of the bottom factor. This is shown in
figure~\ref{Fig:wrapping}(b). We therefore need not make any adjustment, as the
width is enumerated correctly.
Finally, the required `1-unimodal' polygons are obtained by translating the fixed
vertical step below the join to the right edge of the polygon (as shown in red in
figure~\ref{Fig:wrapping}(b)) to ensure that the polygons are self-avoiding.

In conclusion, without needing to make any extra adjustment, wrapping allows us
to enumerate almost-unimodal polygons with a single indent in the left side by
{\em only} enumerating staircase polygons joined to unimodal ones. This is because
the wrapping action also creates polygons that are composed of a unimodal bottom
factor joined to a pyramid top factor.

\section{The 1-convex generating function} \label{s_1-convex}
As an example of the above wrapping technique, we derive the
generating function of 1-convex polygons with their indent on the top arc. It
provides a much simpler derivation of the result than the method used in
\cite{wrgj}.

\begin{figure}
	\begin{center}
	\subfigure[The inclusion case.]
		{\qquad \includegraphics[scale=0.5]{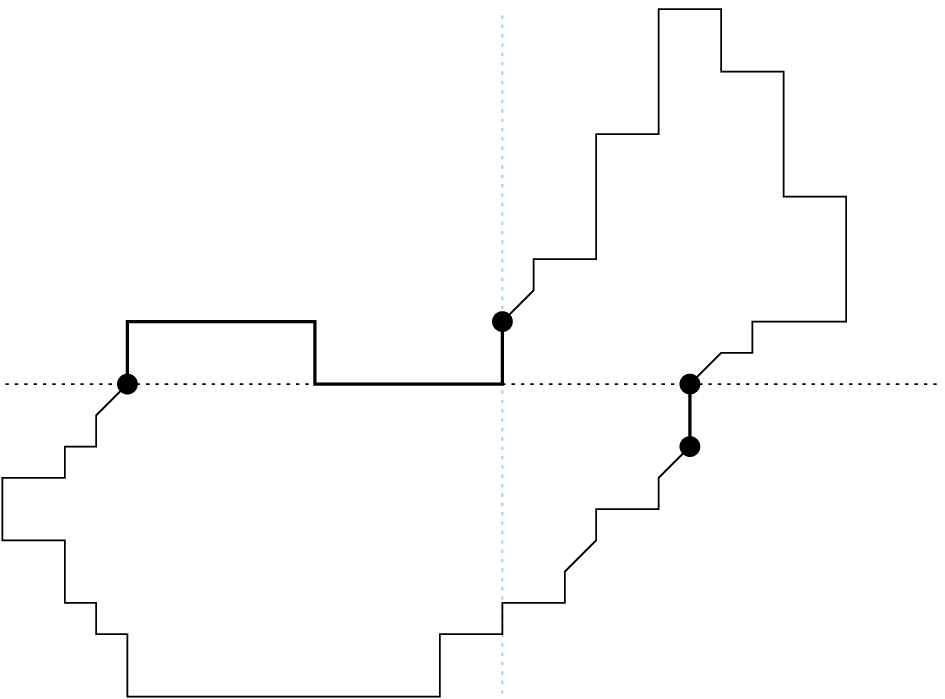} \qquad}
	\subfigure[When the indent extends to the left.]
		{\qquad \includegraphics[scale=0.5]{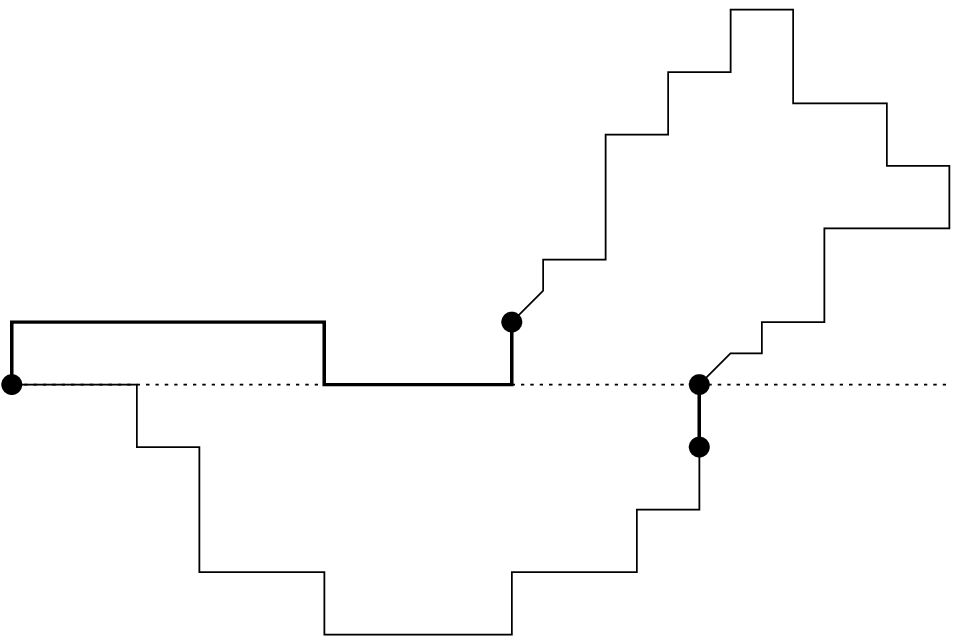} \qquad}
	\end{center}
	\caption{Generalising Lin's factorisation to enumerate 1-convex polygons.}
	\label{Fig:1-convex}
\end{figure}
We begin by adopting Lin's factorisation of 1-convex polygons (see
\cite{Lin:1-convex}) by extending a line along the base of the indent, as in
figure~\ref{Fig:1-convex}. We enumerate these polygons by following the
inclusion-exclusion argument of \cite{BMG:convex}, enumerating all the required
polygons, including those that intersect, and then excluding those that
intersect. The factorisation gives top and bottom unimodal factors that can be
enumerated as walks that may be wrapped such that the bottom factor extends
furthest to the right.

The generating function of polygons of the form shown in part (a) of the figure
can be expressed as
\begin{eqnarray} \nonumber
\fl \sum_{n\geq 1} \frac{1}{x^n} \left(
	\EE{x^{n+1}y^2}{1-x-y}{\left( 1+\frac{y}{1-x} \right)}-u^nv\EE{xy}{1-x-y}{}
\right)  \\ \cdot \left(
	\EE{x^{n+1}y^2}{1-x-y}{\left(\frac{x^2/(1-y)}{1-x^2/(1-y)}\right)^2
	\left(\frac{x}{1-y}\right)^{n-1}}- \frac{u^{n+2}v}{(1-u)^2}\EE{xy}{1-x-y}{}
\right) \nonumber \\
- 2x\SP^3Z^3\left(\frac{1}{1-x}+uZ\right)
  \left(1+\frac{v}{1-u}\right).
\end{eqnarray}
The length of the join is taken to be $n$.
The term in the first line is the generating function for the unimodal top
factor, with a base of at least length $n$. The term $x^{n+1}y^2/(1-x-y)$
enumerates DWs that are folded by the E operator to form the unimodal factor.
We note that due to the vertical symmetry, when this top factor is of height one,
some polygons may be double-counted, giving the term $(1+y/(1-x))$, the
$y/(1-x)$ forcing the top factor in the symmetric case to be of at least height
two.  The second line enumerates the bottom factor and the indent. Again, using
an inclusion-exclusion approach, we fold a walk to form the unimodal polygon and
then exclude the intersecting cases. As the width along the join is already
counted in the top factor, we adjust for the length of the join with the term
$(x/(1-y))^{n-1}$. We then fix the horizontal steps of the indent factor, which
cannot be folded, giving the term
\[
	\left(\frac{x^2/(1-y)}{1-x^2/(1-y)}\right)^2,
\]
rather than the expected $x^2/(1-x)^2$, which can be folded.
We recognise the last term in both of the first two lines as the exclusion cases
when the polygons intersects in the top-right or bottom-left corner. The last
term comes from the possibility of intersection in the bottom-right corner.

From the wrapping principle outlined in the previous section, when the top
factor has the form of a pyramid, the fixed steps along the join may have been
folded. This then wraps the bottom factor, making it convex in shape, such that
it extends further to the right than the top factor.

We finish by expanding the brackets, so that we may evaluate the sums and write
the expression in a closed form.
However, the indent may extend further to the left than the bottom factor, as
shown in part (b) of the figure, and therefore requires an adjustment term.
And so, moving all terms incorporating $n$
into the $E$ operators and expanding the brackets in the summand, we can then
complete the summation. This then allows us to add the term $x^*/(1-x^*)$ that
enumerates the indent for the required adjustment to give the following expression
for the generating function:

\begin{eqnarray} 
\fl \EE{x(1-x)y^*}{(1-x)^2-y^*} {
	\lr{1+\frac{y^*}{(1-x)^2}}\lr{\frac{x^*}{1-x^*-y}}^2
	\lr{\frac{y^2(1-y)^2}{(1-y)^2-x^*}+\frac{x^*y^2}{1-x^*}} \frac{x}{1-x-y}} \nonumber \\  
- \frac{4xyv}{\Delta} \EE{}{}{ \lr{\frac{x^*y}{1-x^*-y}}^2
	\lr{\frac{(1-y)^2}{(1-y)^2-x^*}+\frac{x^*}{1-x^*}} \frac{u^*}{1-u^*-y}} \nonumber \\  
 - \frac{2xyu^2v}{(1-u)^2\Delta} \EE{x(1-x)y^*}{(1-x)^2-y^*}{
	\lr{1+\frac{y^*}{(1-x)^2}} \frac{x}{1-x-v^*}} \nonumber\\ 
+ 2 v \SP Z \lr{ \frac{2x\SP}{\Delta}}^2
- 2x\SP^3Z^3\left(\frac{1}{1-x}+uZ\right)
  \left(1+\frac{v}{1-u}\right).
\end{eqnarray}

\section{Derivation of the 2-convex generating function}

We generalise Lin's factorisation of 1-convex polygons by extending a factorisation
line along the base of \emph{each} indent. This allows each case to be enumerated by
joining factors along these lines. When the indents are 
in the same direction, there are three main factors plus two indents. When the 
indents are in different directions, we divide the lattice into four quadrants. 
We say the quadrant in the top-right is the first, and order the remaining
quadrants in an anti-clockwise fashion.

To break-up the problem into enumerable parts, similar to Lin, we adopt a
`divide and conquer' approach and classify  sub-classes of 2-convex
polygons according to the relative direction and position of the indents (that
is, which side they lie on). We obtain the generating functions of
symmetric classes by reflection and rotation.  Without loss of generality 
we assume that one indent is vertical and on the top-left side. This leaves us
with one of nine cases: Firstly the indent can have depth two, secondly when there
are two indents each of depth one, we have to consider eight 
combinations of the direction and location of the second indent.  
However, the two cases where the second indent is in a different
direction and on one of the adjacent sides are equivalent after a rotation.
We are therefore left with eight distinct cases to evaluate.

There is a possible ambiguity when an indent factor is adjacent to the MBR. For
example, in the case where there is a vertical indent on the base of the polygon
as well as on the top-left side, the indent on the base may be considered as
on either side of the bottom arc. We arbitrarily chose that such cases be
enumerated by the class whose indented sides are closest together. This example
is therefore enumerated by the case where the indents are on adjacent sides,
and not on opposite sides of the polygon. We now briefly consider these eight
cases in turn.

\subsection{Case 1: a single 2-deep indent}
Almost-convex polygons with a single indent we refer to as `bimodal' due to
the two modes (in the same sense as `unimodal') adjacent to the indent. In section
\ref{s_1-convex} we enumerated 1-convex polygons, which are the simplest case of
bimodal polygons. Bimodal 2-convex polygons, which have a single 2-deep indent,
can be enumerated {\em mutatis mutandis}. For the generating function, see
Section~2 of \cite{JJG07}.

\subsection{Case 2: indents in the same direction on the same side}

\begin{figure}
	\begin{center}
		\includegraphics[scale=0.5]{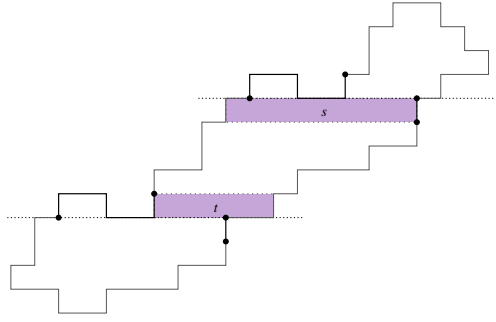}
	\end{center}
	\caption{The form of a 2-convex polygon with both indents on the same side.}
	\label{Fig:2-convex}
\end{figure}

When the indents are on the same side (on the top-left), we join unimodal top
and bottom factors to a staircase factor in the middle. This implies that the top
factor extends furthest to the right. (See figure~\ref{Fig:2-convex}.) This
creates three blocks separated by two factorisation lines. Each pair of blocks
are joined by matching the top and base edges along the factorisation line dividing
them. Since there is an indent factor in the top block, horizontal steps 
must be added on either side of the indent so
their perimeters match. In this way the length of the join,
counted by the parameter $s$, is equal to the top perimeter of
the middle staircase factor. Whenever the join is defined in this
way,  in order to enumerate the indent, a term $s^2/(1-s)^3$ is included on one 
side of the Hadamard product in the expression for the generating function. Indeed,
this term appears in the
majority of expressions involving 2-convex polygons. Fortunately, such joins
can be re-expressed in terms of the first three moments of the generating
function with the indent omitted.

In order to match lengths correctly
along the bottom join (defined as the length of the base of the middle
staircase factor and enumerated by the parameter $t$) we append horizontal
steps to the bottom factor.  In this example, the length of both joins are
defined as the length of the middle factor along the
factorisation lines, as indicated by the shaded regions in figure~\ref{Fig:2-convex}. 
The choice for the definition of the join is determined by
the generating functions of the top and bottom factors, which are simple
rational expressions in terms of $s$ and $t$. This is
generally simpler than trying to adjust the middle factor to match the other factors.

In using standard inclusion-exclusion techniques for enumerating the top and bottom 
unimodal factors, wrapping ensures that the cases where the middle or bottom factor 
extends furthest to the right are included. This means that the bottom 
factor may intersect (forming a unimodal loop) to the bottom-right. This is a powerful 
example of the robustness of the wrapping methodology -- we fold the top factor, wrapping not
just the middle factor, but the bottom factor as well, such that it can be the
one to extend furthest to the right.

We now complete the enumeration to serve as an example for the following cases.
We only give some details here, as there are dozens of very complicated formulae
in all, each using the same principles in their derivation. We outline each of
the remaining cases to the minimum extent that would be necessary to reproduce
the results.
We break the problem into two parts: when the indents are at the same height
and when they are not. The latter case is enumerated by the following
expression, and the former follows \emph{mutatis mutandis}.

\begin{eqnarray}
\fl	\lr{ \EE{sxy^3}{(x-s)(1-x)(1-x-y)}{}
		- \frac{2xy}\Delta \cdot \frac{sv}{1-s-v} } \frac{s^2}{(1-s)^3}
	\ \odot_s\ \bar S(s,t)\ \nonumber \\
	\odot_t\ \frac 1{1-t} \EE{ty^2}{1-t-y}
	{\lr{\frac{(1-y)^2}{(1-y)^2-x^*}+\frac{x^*}{1-x^*}} \lr{\frac{x^*}{1-x^*-y}}^2\,}
	\\ \nonumber
	- \frac{2xy}\Delta \lr{
		\lr{ \frac{y^4}2 \dd y \T^2 \dd y \frac 1{y^2} \lr{xyZ - \frac{xy}{1-x}} }
		+ v^2 \SP Z^2 \T^2 \lr{\SP Z + \frac v{1-x}}
	},
\end{eqnarray}
where $\bar S$ is the  generating function of staircase polygons by base and top
perimeter.

The expression to the left of the join in the first line enumerates the top
factor, with the indent enumerated by $s^2/(1-s)^3$. The first term in the $E$
operator enumerates the (possibly intersecting) unimodal folded walks, with the
term $s/(x-s)$ counting the fixed steps along the join (weighted by $xs$, but
divided by $x^2$ to adjust for the width already enumerated by the middle
factor). The term $y/(1-x)$ ensures that it is at least of height one. The
inclusion-exclusion principle then lends us to exclude the intersecting case
enumerated in the second term.
The second line enumerates the bottom factor, obtained by folding a
stack polygon (a reflection in the horizontal axis). This ensures that
the indent is not wrapped. The first term is the part of the stack polygon under
the join. The last term is the part under the indents. The second term
enumerates the part of the polygon that extends to the left, including the
possibility that the bottom factor is a pyramid and that the indent factor
extends furthest to the left.
Finally, the terms in the last line enumerate the exclusion cases.
The first term counts the polygons that intersect in the
bottom-left corner. The 2-unimodal factor is enumerated by distinguishing the
heights where the indents are placed, thus explaining the derivatives. The center
term of the derivative is the generating function for
unimodal polygons of at least height two. The last term counts polygons
intersecting in the bottom-right corner. This term is obtained by summing the
generating functions for each possible configuration of the indents, using
the known generating function for staircase polygons with fixed steps in the corner.

\subsection{Case 3: indents in the same direction on adjacent sides of the same arc}

\begin{figure}
	\begin{center}
	\subfigure[The left-side indent is above the corner indent.]{
		\includegraphics[scale=0.38]{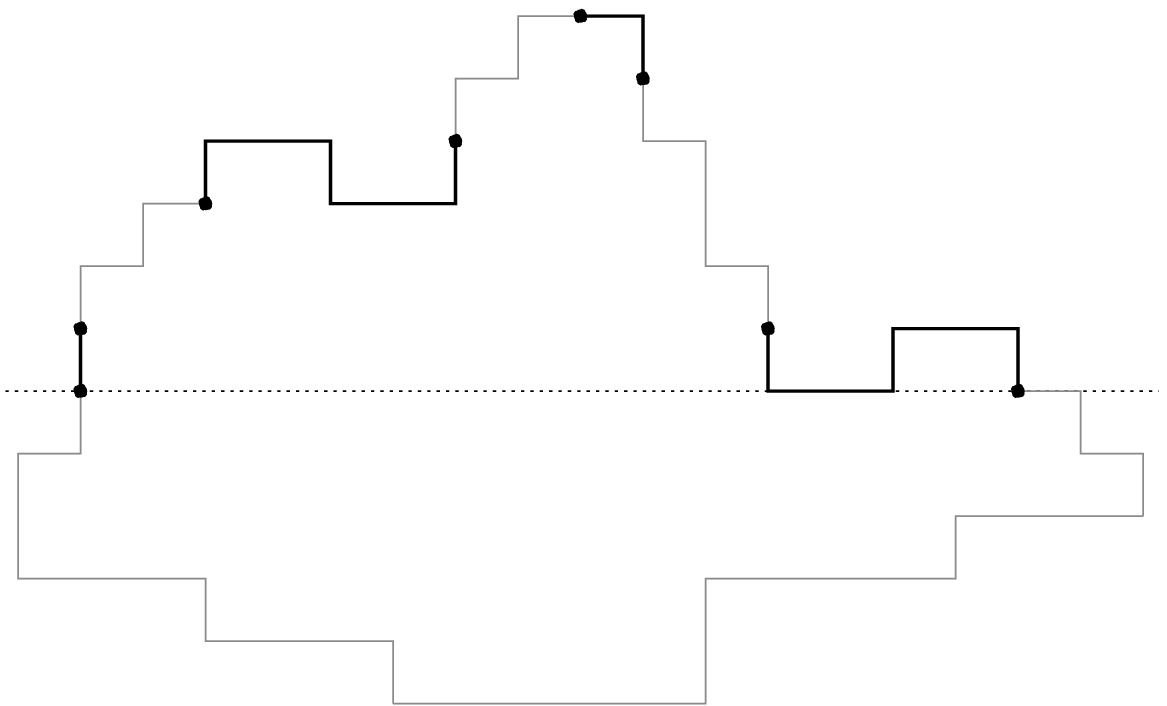}} \quad
	\subfigure[The left indent is below the corner indent.]{
		\includegraphics[scale=0.38]{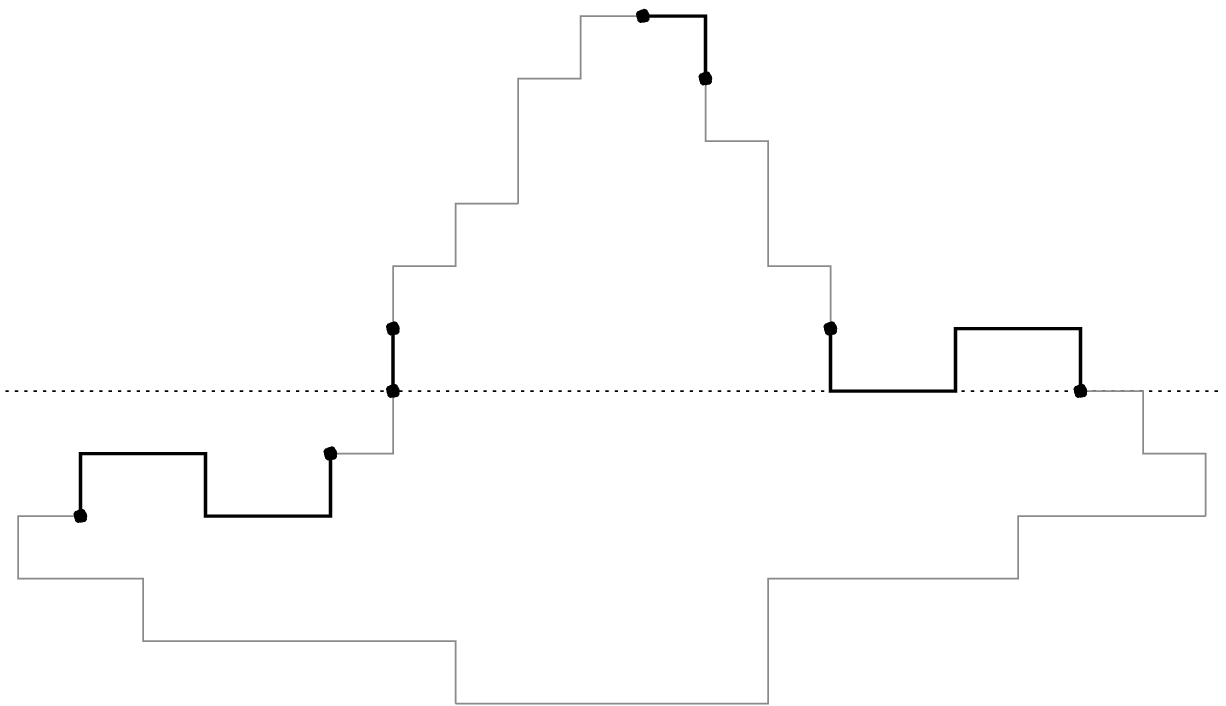}}  \quad
	\subfigure[Both indents are on the same side.]{
		\includegraphics[scale=0.38]{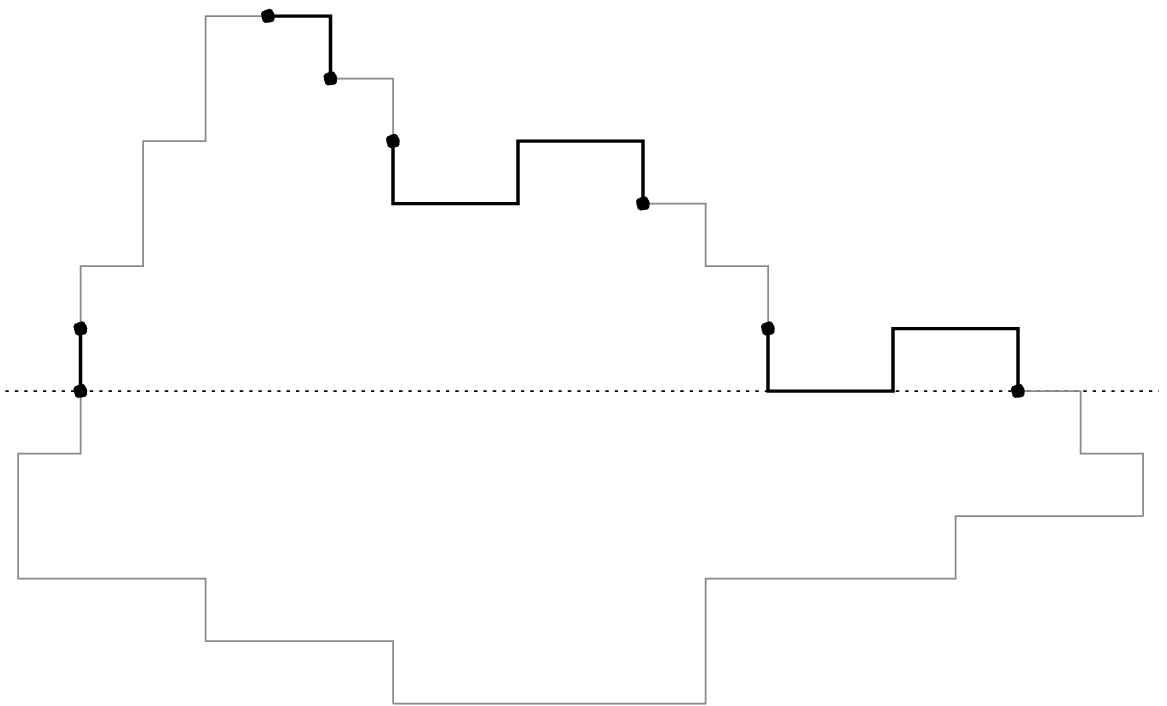}}
	\end{center}
	\caption{The form of 2-convex polygons with two vertical indents on adjacent
			 sides of the top arc.}
	\label{Fig:case3}
\end{figure}

In this case the top factor must be a pyramid.
We break the enumeration into two parts depending on whether both indents are at the same
height or at different heights.
For the latter case,  illustrated in figure~\ref{Fig:case3}, we assume
that the top-left indent is higher than the top-right one. The symmetrical
case can be obtained by reflection.

In part (c) of the figure, we reproduce the situation from Case 2 (the
indents are on the same side) where the top factor is a pyramid (such that the
middle factor is wrapped). Reflecting the indented top pyramid factor
vertically we recover the form of the polygons shown in part (a).
Alternatively, drawing the factorisation line under the top indent, flipping
both the indent and the top factor gives the polygons shown in part (b).

We can therefore use the expression for the generating function for Case 2, with
a slight modification to the first term, which becomes
\begin{equation}
	\EE{sy^2}{1-x-y}{\lr{\frac x{x-s} - \frac{1-y}{1-y-s}}} +\frac{P(s,y)}{1-x},
\end{equation}
where $P(x,y)$ is the pyramid generating function.
The pyramid term comes from the possibility that the indent 
extends furthest to the right. The last term in the $E$ operator excludes those
cases from Case 2 that were not wrapped and therefore did not have the pyramid top
factor now required.

\subsection{Case 4: indents in the same direction on opposite sides}
\begin{figure}
	\begin{center}
	\begin{tabular}{ccc}
		\subfigure[The top indent is below the bottom indent.]
		{\raisebox{4pt}
		{\includegraphics[scale=0.4]{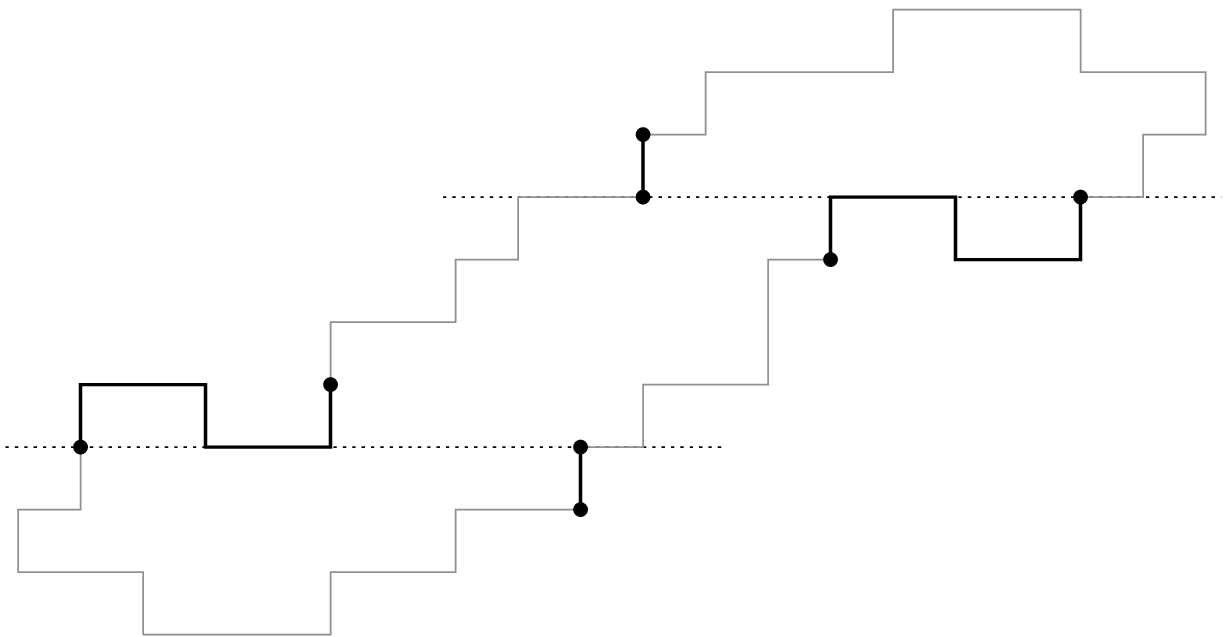}}} &
		\subfigure[The top indent is level with the bottom indent.]
		{\raisebox{16pt}
		{\includegraphics[scale=0.45]{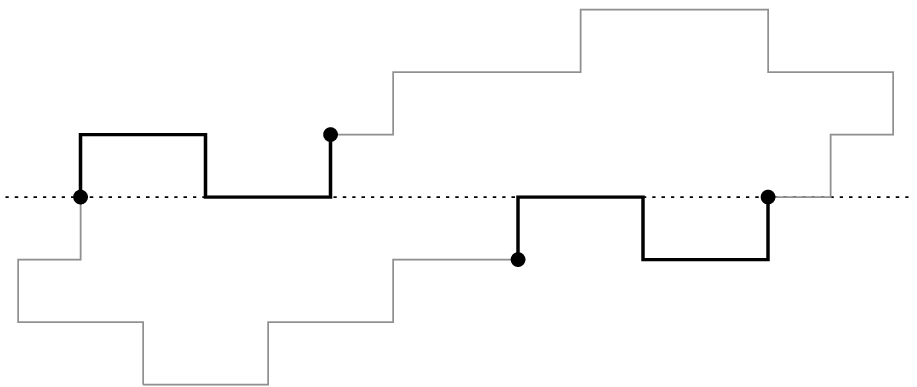}}} &
		\subfigure[The level indents are inter-weaved.]
		{\raisebox{16pt}
		{\includegraphics[scale=0.45]{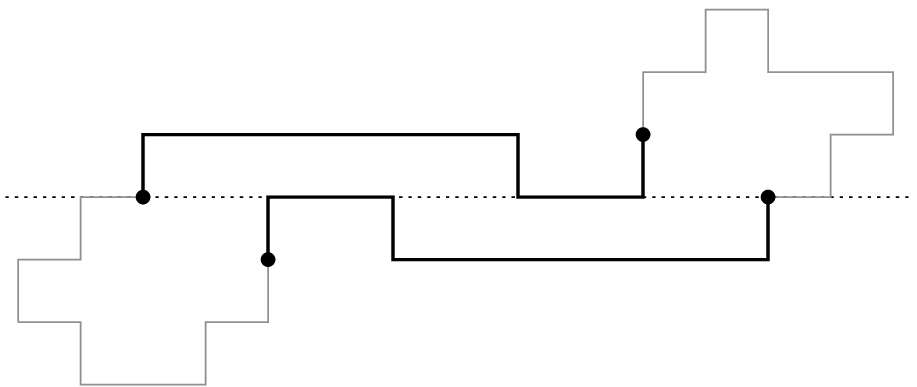}}} \\
		\subfigure[The top indent is above the bottom indent.]
		{\raisebox{4pt}
		{\includegraphics[scale=0.4]{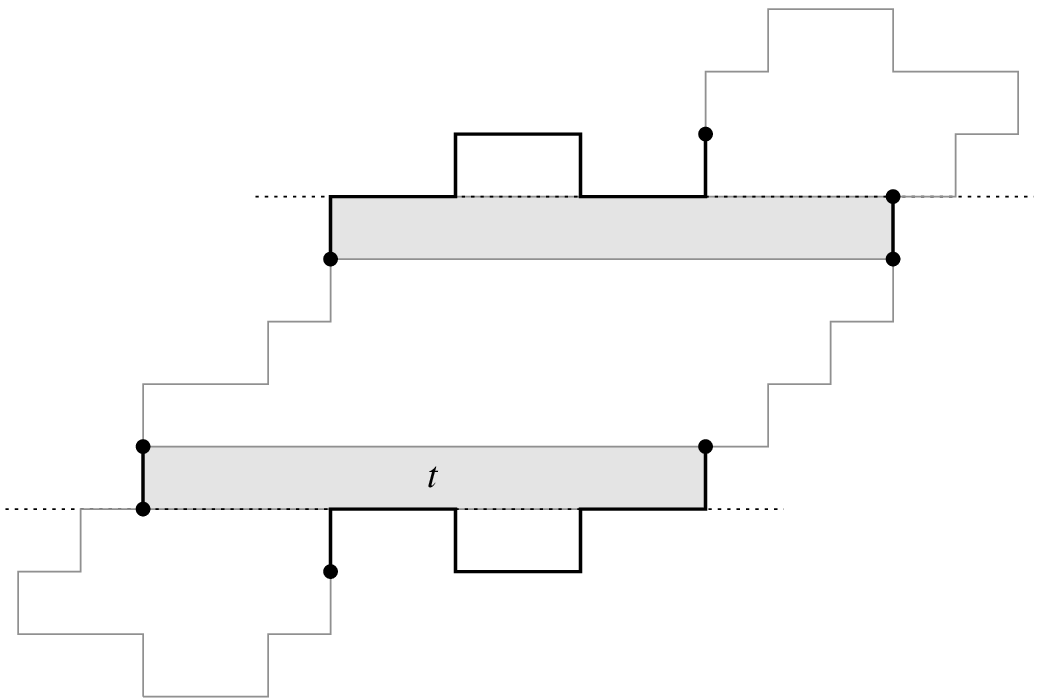}}} &
		\subfigure[The middle factor is further to the right than the top factor.]
		{\raisebox{25pt}
		{\includegraphics[scale=0.4]{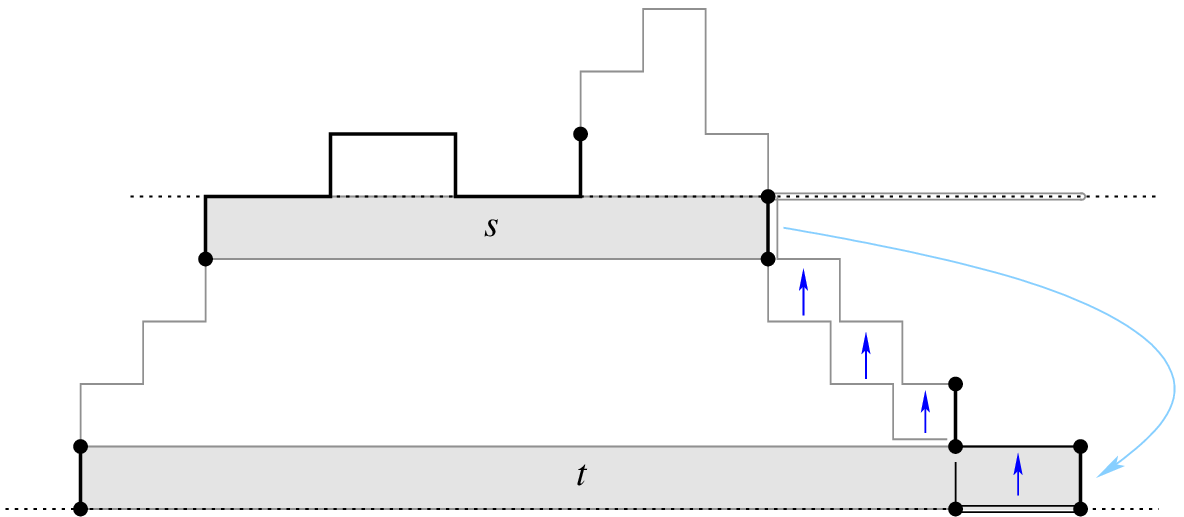}}} &
		\subfigure[The indent factor is farthest to the right.]
		{\raisebox{4pt}
		{\includegraphics[scale=0.4]{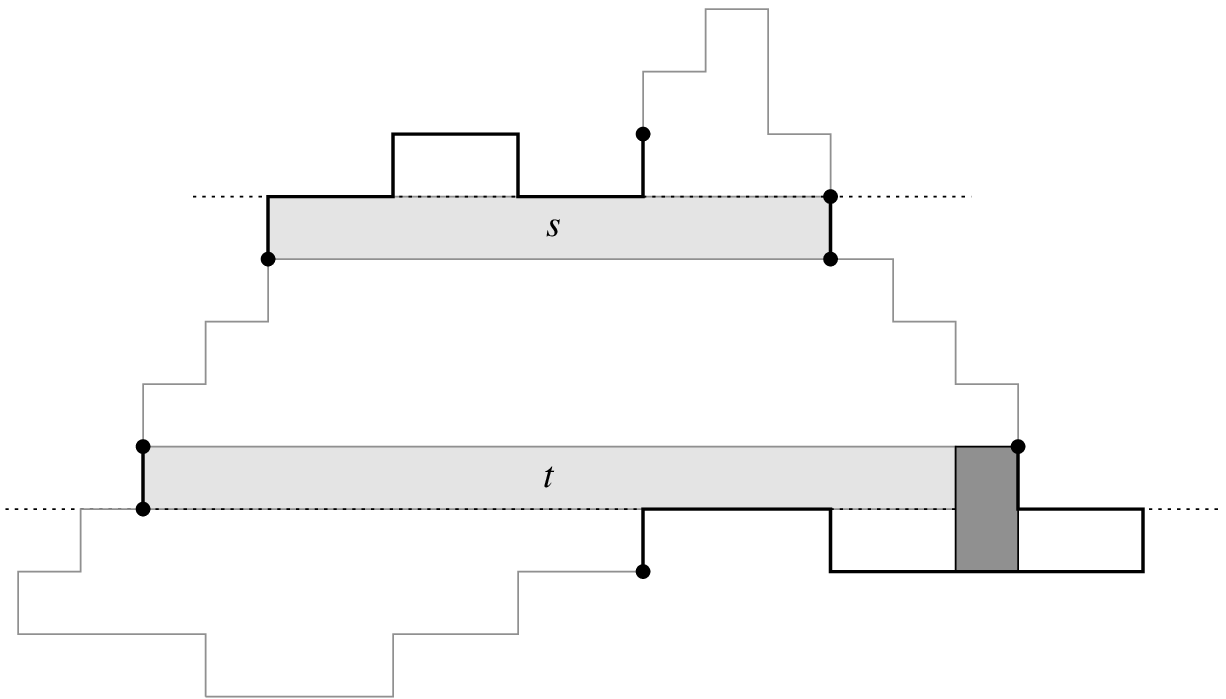}}}
	\end{tabular}
	\end{center}
	\caption{The form of 2-convex polygons with indents on opposite sides.}
	\label{Fig:case4}
\end{figure}

We divide this case into three parts according to whether the top-left indent is
above, next to or below the bottom-right indent (see figure~\ref{Fig:case4}.)
When the top indent is above the bottom one, and either the top or
bottom factor is of height one, the indents can be considered to be on
adjacent sides, and we do not include them in this case. Therefore, all of the
above calculations must be done for top and bottom factors of height at least
two.

When the top indent is below the other (see part (a) of the figure), the
enumeration does not require wrapping and has no complications; it can be
obtained by simply joining the three appropriate factors together.
When the indents are next to one another (see part (b) of the figure), the
enumeration is even simpler. However, when the indents are level
the top indent may be either to the left \emph{or the right} of the bottom-right indent.
The latter case means that the indents interweave (see part (c) of
the figure).

Now consider the situation where the top-left indent is above the other one.
When wrapping, in order to keep the polygon self-avoiding, we usually
translate the fixed vertical step from below the join to the fold (see
figure~\ref{Fig:wrapping}(b)). However, when wrapping the bottom factor as well
we need to make an adjustment to the length of the bottom join (see
part (e) of the figure). Moreover, when the middle factor is a pyramid, it is possible 
that the indent extends furthest to the right (see part (f) of the figure).

This case presents a complication not met earlier in our enumeration. Proceeding
as usual would involve simultaneously joining an indent and unimodal factor to both
the top and bottom of a staircase middle factor. When the middle
factor extends furthest either to the left or right  wrapping 
generates these polygons. However, this requires wrapping both sides of the middle
staircase factor independently, which is not possible as the calculations
diverge. (For further discussion, we refer to Section 7.2.1 of \cite{JJG07}.)
We can therefore only use wrapping on one side, and must break up the
calculation of the bottom factor into two parts depending on whether it is unimodal
or pyramid.

Finally, note that when the middle factor extends furthest to both the left
and right, it is convex in form. This requires
the enumeration of convex polygons according to both  base and top-most
horizontal segment.
This can be achieved by adopting a `divide and conquer' approach
and joining pyramid factors to unimodal ones, or by solving recurrence
relations. (The generating function by perimeter and area, solved by the
aforesaid recurrence relations, already exists in the literature \cite{Bousquet:96}.)

\subsection{Case 5: indents in the same direction on adjacent sides of opposite arcs}
\begin{figure}
	\begin{center}
		\includegraphics[scale=0.5]{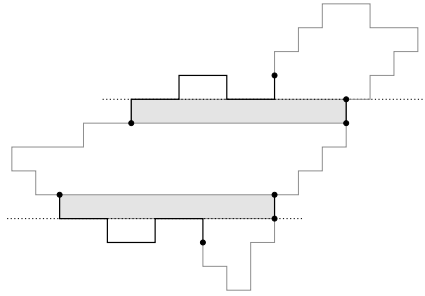}
	\end{center}
	\caption{The form of 2-convex polygons with indents in the same direction on
	adjacent sides of opposite arcs.}
	\label{Fig:case5}
\end{figure}

This case is similar to the previous one, when the bottom factor was a pyramid.
Reflecting the bottom factor and indent  the
correct form is produced as shown in figure~\ref{Fig:case5}.
We can proceed by joining a unimodal top factor to a unimodal middle factor,
which in turn is joined to a pyramid bottom factor. Wrapping again generates the
cases where the bottom factor extends furthest to the right.
The enumeration of the different parts of this case follows the same approach as
the previous case, except that the top and bottom factors may be of
height one. Also, the interwoven case is somewhat different, as the
top and bottom factors are joined directly.

In defining these eight cases we have assumed that the top
indent is on the top-left side. When multiplying by two to obtain the generating
function for the symmetric case (where the top indent is on the top-right side),
we double-count polygons whose top and bottom indents are adjacent to the MBR
(i.e. they form the topmost segment and the base).  We must therefore adjust for
this case when adding up the generating functions at the end.

\subsection{Case 6: indents in different directions on the same side}
\begin{figure}
	\begin{center}
	\subfigure[The indents form a concave region in the corner.]
		{\includegraphics[scale=0.4]{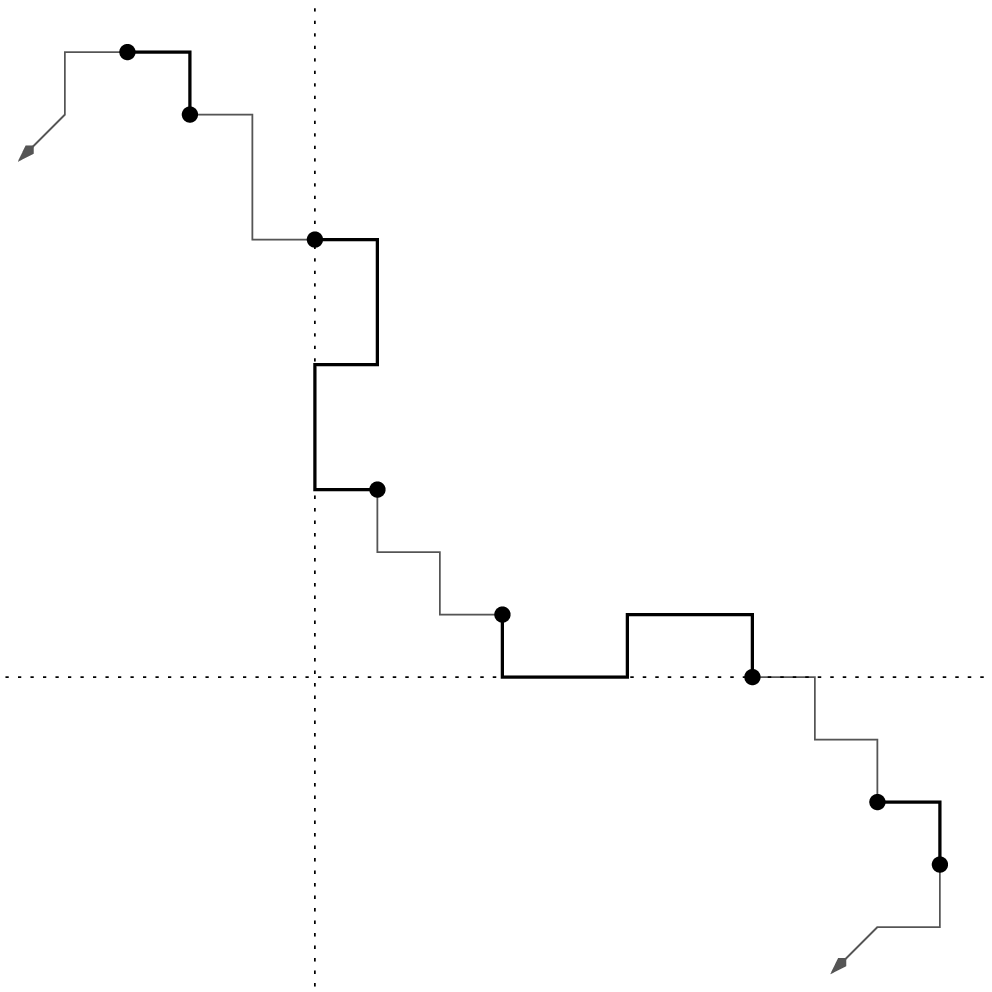}
		 \raisebox{12pt}{\includegraphics[scale=0.4]{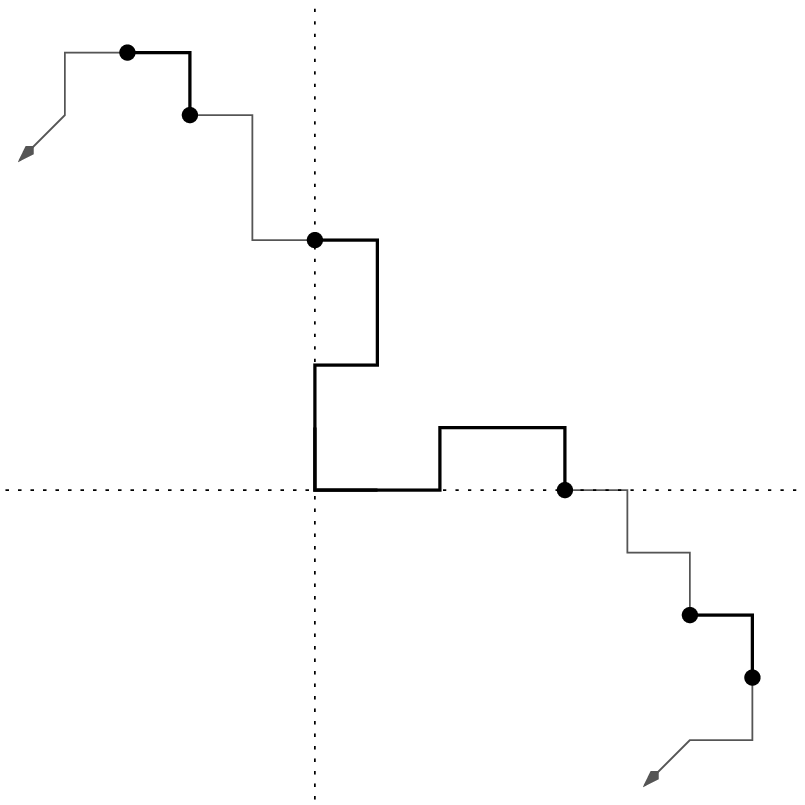}}} \quad
	\subfigure[The indents form a convex region in the corner.]
		{\includegraphics[scale=0.4]{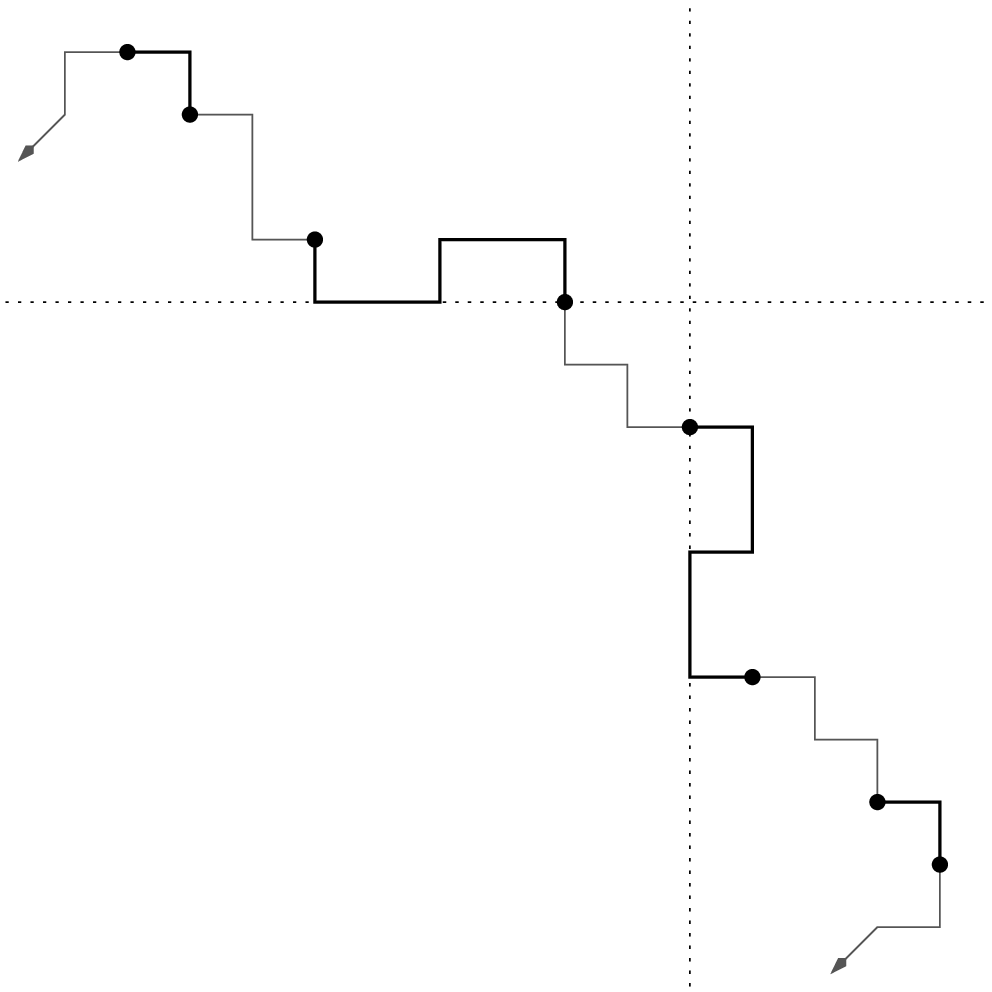}
		 \raisebox{12pt}{\includegraphics[scale=0.4]{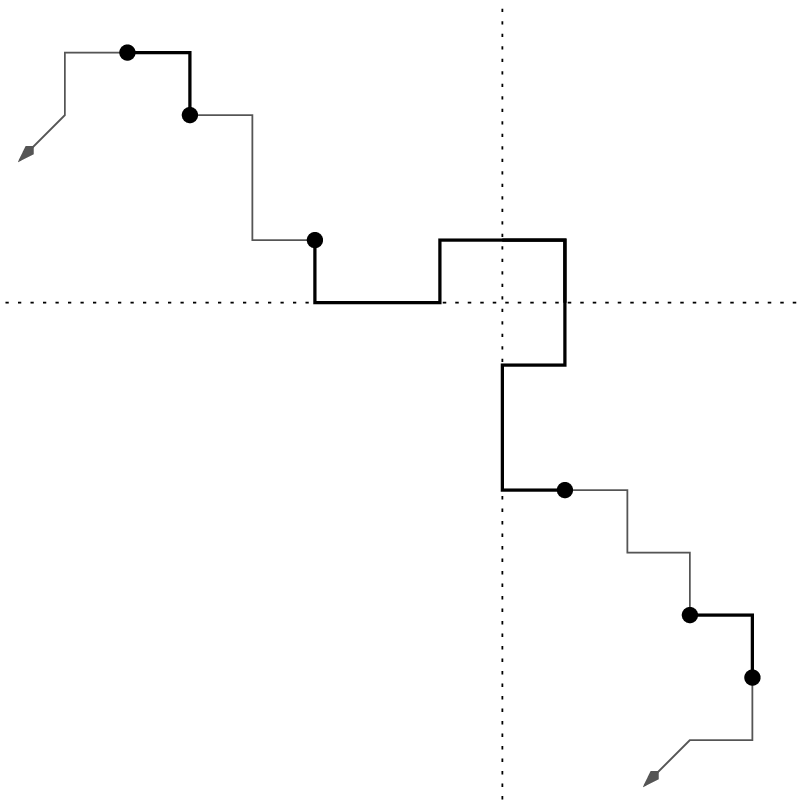}}}
	\end{center}
	\caption{The form of 2-convex polygons with indents on the same side in
			 different directions.}
	\label{Fig:case6}
\end{figure}
When the indents are in different directions and on the same side they
form either a locally convex or concave region, depending on their
order. This is shown in figure~\ref{Fig:case6}.

First consider the locally concave case.
Aside from the walk between the two indents, the polygon can either enter the first or
the third quadrant, but not both. 
If the polygon does not enter the third quadrant then it has only three factors:
two unimodal ones joined to the staircase factor with the indents.
When it passes through the third quadrant, 
the indents and the connecting walk of the first quadrant are joined to unimodal factors in
the second and fourth quadrants. These are in turn joined by a directed walk in
the third quadrant. Making the expressions for the second and fourth quadrants
factor into two parts, one independent of the horizontal join, the other
independent of the vertical join, allows the calculation of the generating
function. This is a good example of how wrapping works in both directions to
generate all the required polygons.

When the indents form a convex region, the polygon does not enter the first
quadrant, except in the case where the indents intersect and there are two steps in
the first quadrant (as shown in the second diagram in part (b) of the figure). The
two cases are evaluated separately, but the second is a simple version of the
first. A unimodal factor and an indent in each of the second and
fourth quadrants are joined to a staircase factor in the third quadrant. The evaluation is
straightforward and wrapping gives the remaining polygons.

\subsection{Case 7: indents in different directions on adjacent sides}
  
\begin{figure}
	\begin{center}
	\subfigure[The top indent is above the one in the corner.]
		{\includegraphics[scale=0.4]{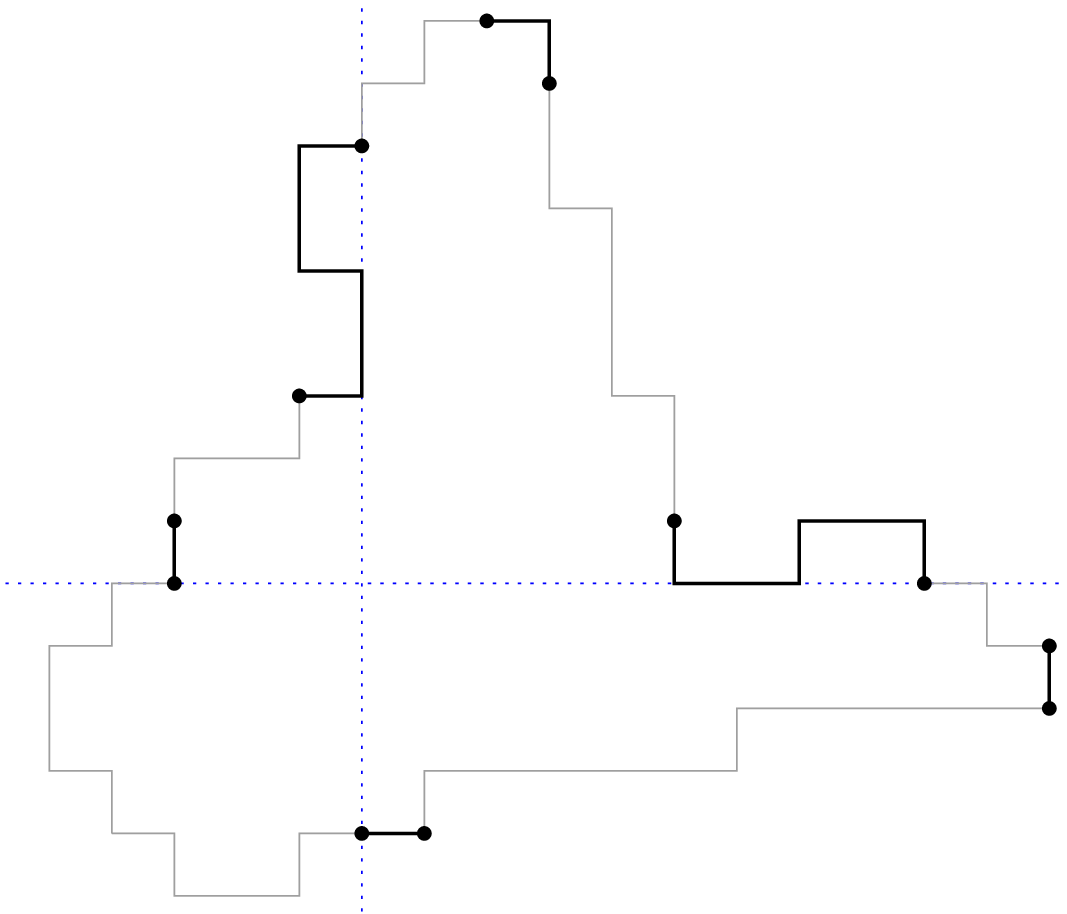}} \qquad
	\subfigure[The left indent is next to the vertical one.]
		{\includegraphics[scale=0.4]{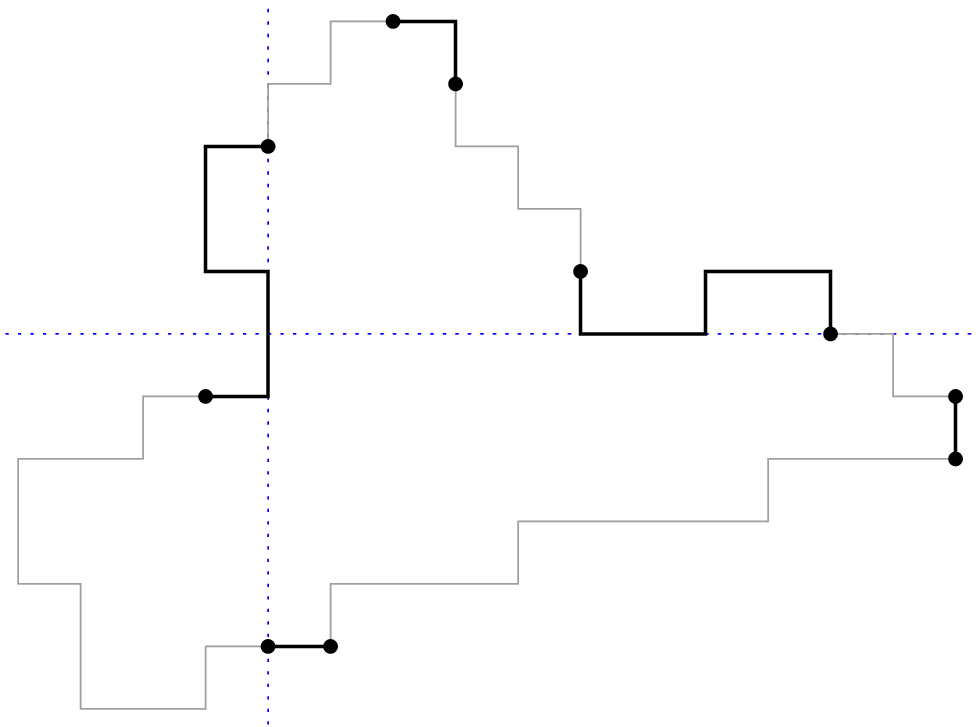}} \qquad
	\subfigure[The left indent is below the other.]
		{\includegraphics[scale=0.4]{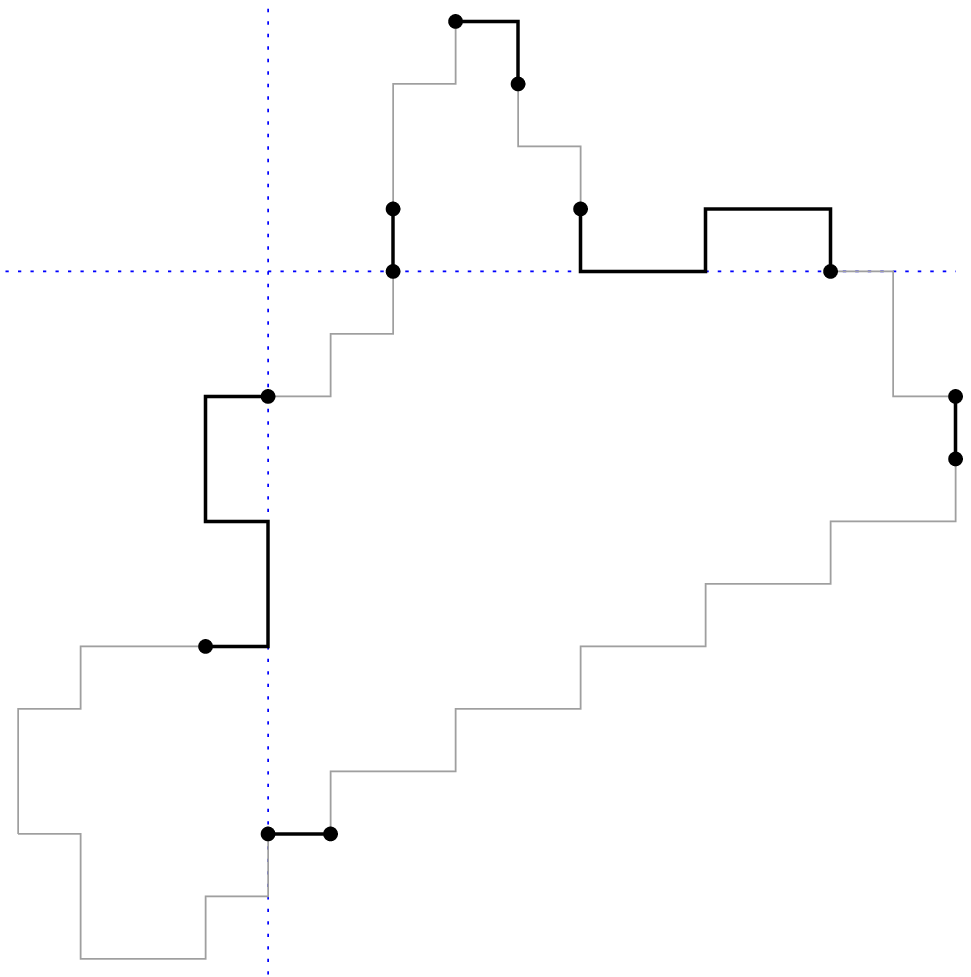}}
	\subfigure[The indent factor is adjacent to the vertical indent.]
		{\includegraphics[scale=0.4]{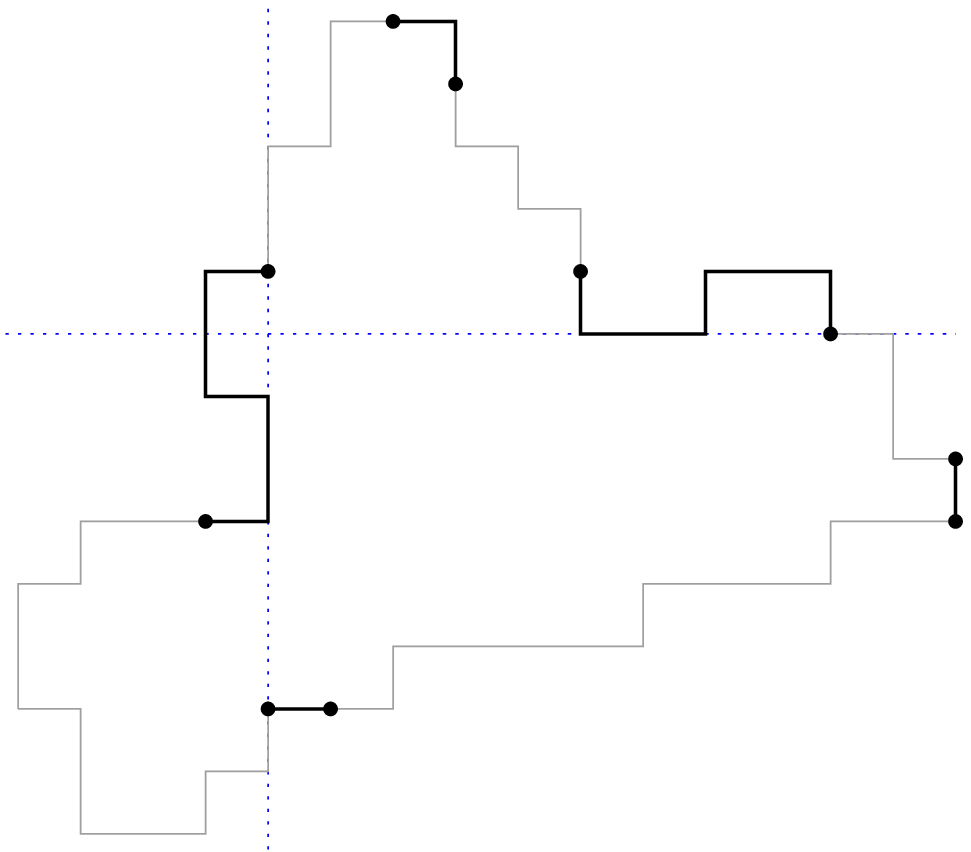}
		 \quad
		 \includegraphics[scale=0.4]{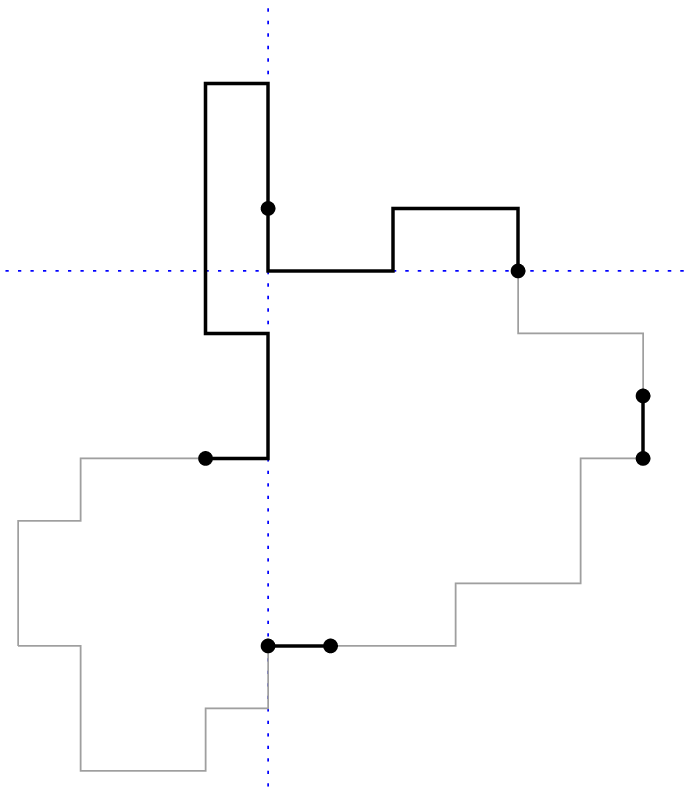}} \qquad
	\subfigure[The polygon does not enter the third quadrant.]
		{\includegraphics[scale=0.4]{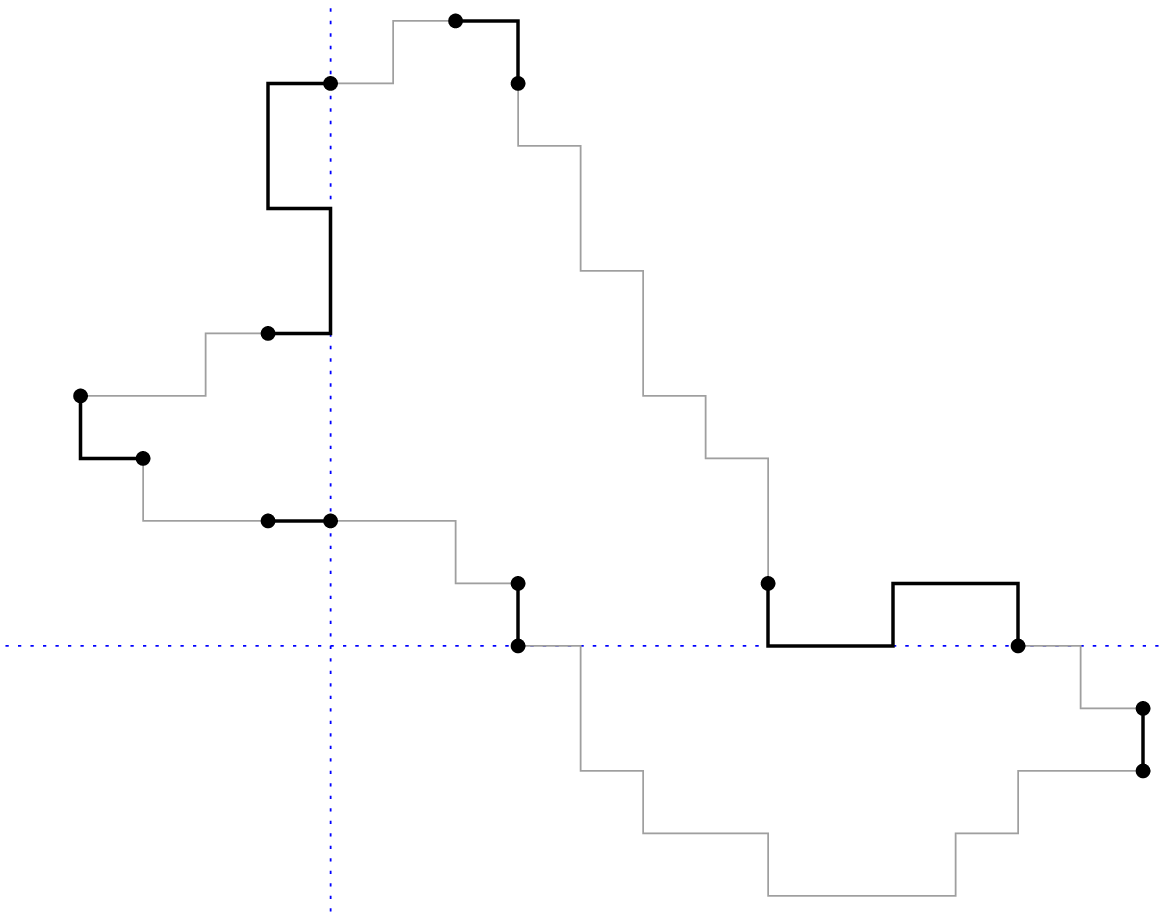}}
	\end{center}
	\caption{The form of 2-convex polygons with indents in different directions
	on adjacent sides.}
	\label{Fig:case7}
\end{figure}

We enumerate this case by breaking it up into parts classified by the
relative height of the indents. We say that the indents are {\em next to} each
other if the {\em holes} formed overlap in height. If the vertical projection of
the {\em humps} formed
by the indent overlap, we say that they are {\em adjacent}.  And so, the horizontal
indent is either 'above', 'next to', 'adjacent to' or 'below' the vertical indent
(see figure~\ref{Fig:case7}).
The cases where the polygon extends furthest to the left in the third
quadrant are all evaluated using the usual inclusion-exclusion and wrapping
arguments.

When the horizontal indent is above the vertical one, the polygon may
extend furthest to the left in the second quadrant (as shown in part (e) of the
figure).
We  need to enumerate 1-unimodal polygons whose indent is in
the corner according to their base. This can be done by joining an almost-pyramid
polygon to a unimodal one. These are then joined to the bottom factor together with an
indent factor,
noting that the indent may extend further to the right than the bottom factor.

\subsection{Case 8: indents in different directions on opposite sides}

\begin{figure}
	\begin{center}
	\subfigure[The top indent is above the one in the corner.]
		{\includegraphics[scale=0.4]{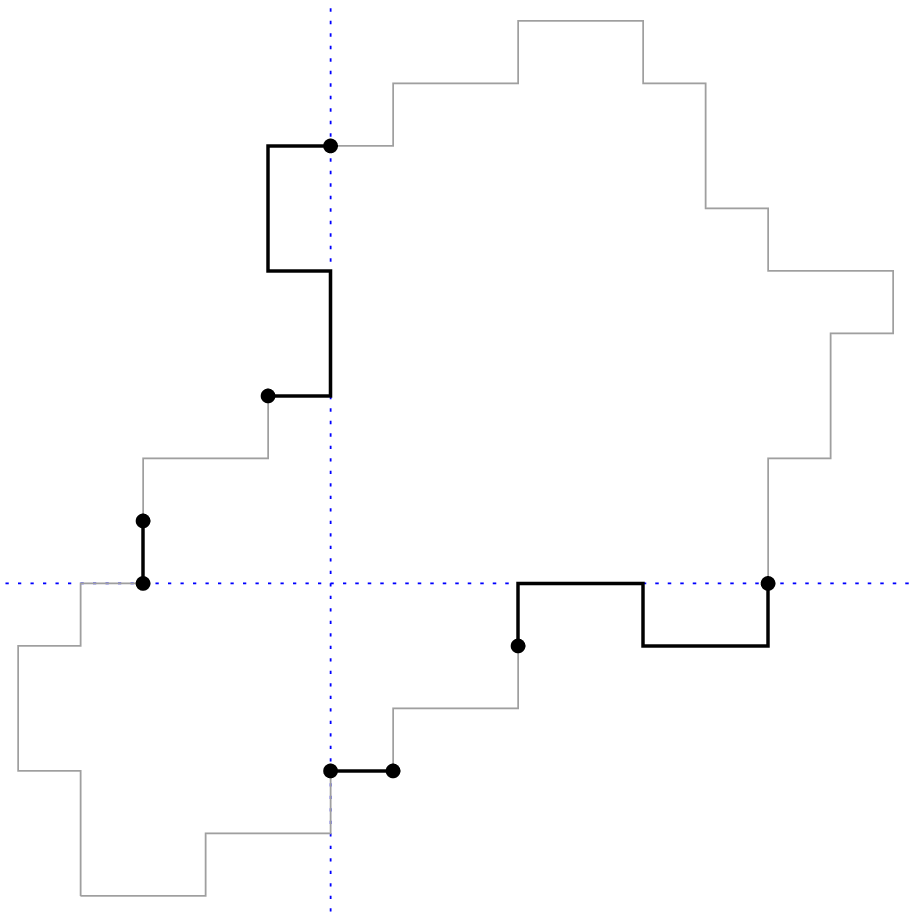} \quad
		 \includegraphics[scale=0.4]{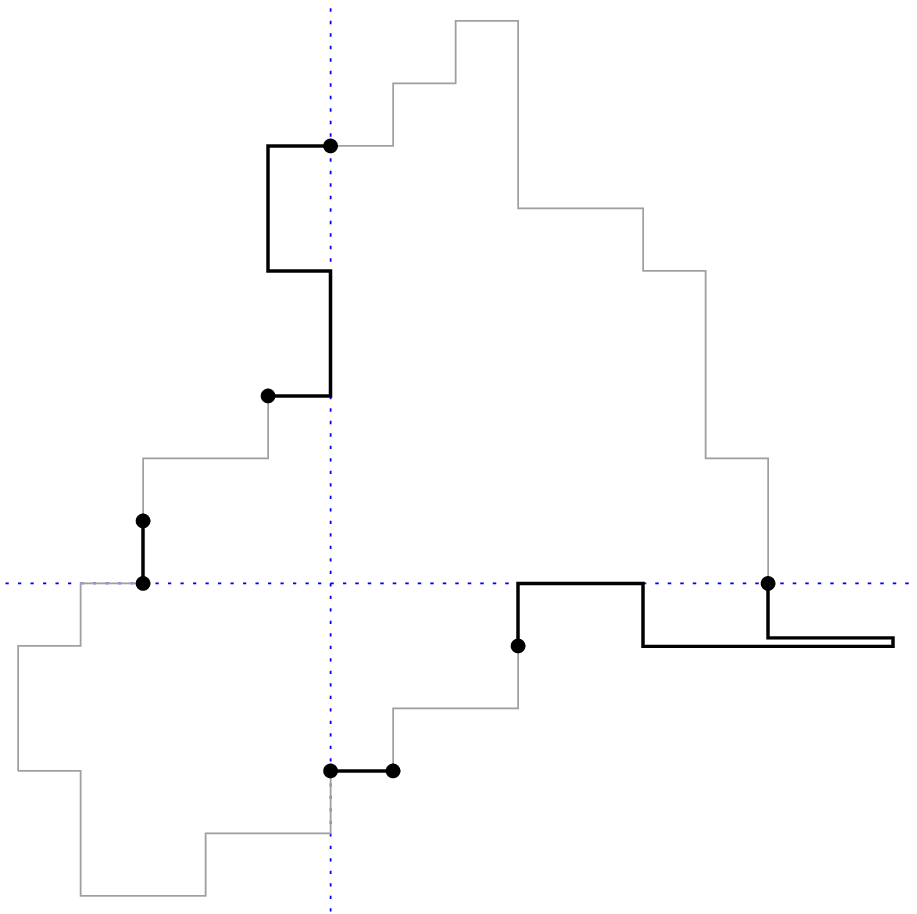}} \qquad
	\subfigure[The left indent is next to the vertical one.]
		{\quad\includegraphics[scale=0.4]{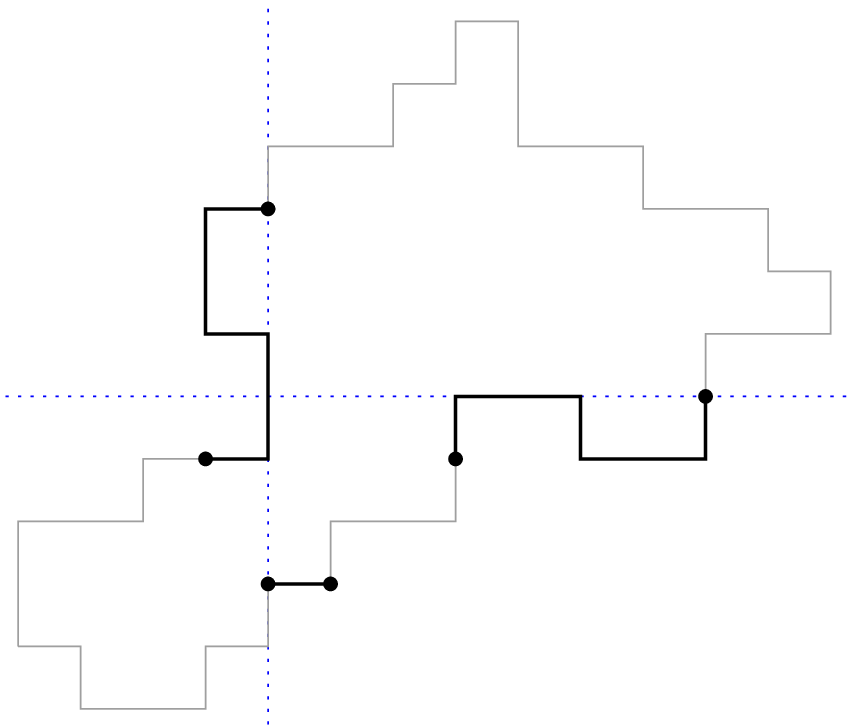}\quad}
	\subfigure[The indent factor is adjacent to the vertical indent.]
		{\quad\includegraphics[scale=0.4]{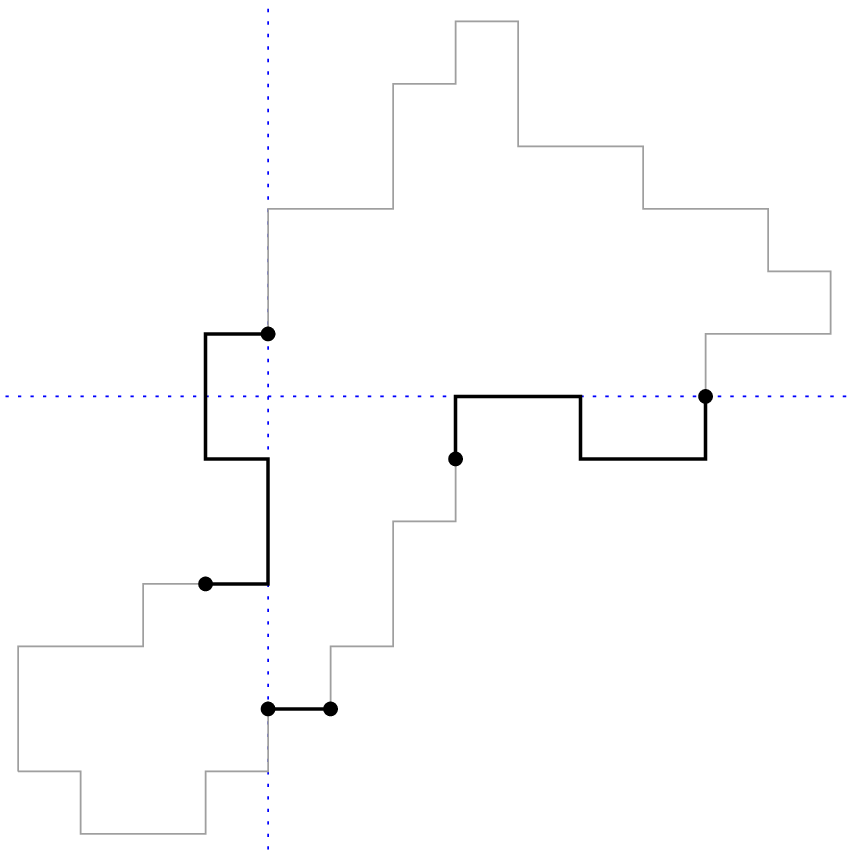}\quad} \quad
	\subfigure[The left indent is below the other.]
		{\includegraphics[scale=0.4]{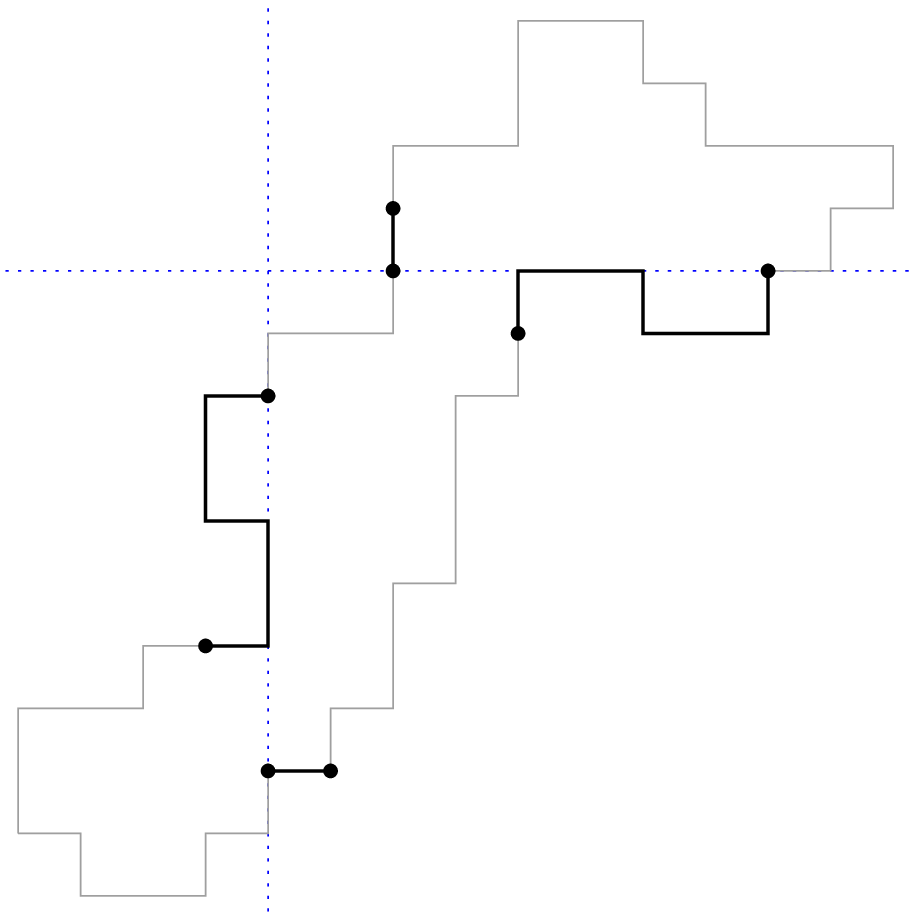}} \qquad
	\subfigure[The polygon does not enter the third quadrant.]
		{\includegraphics[scale=0.4]{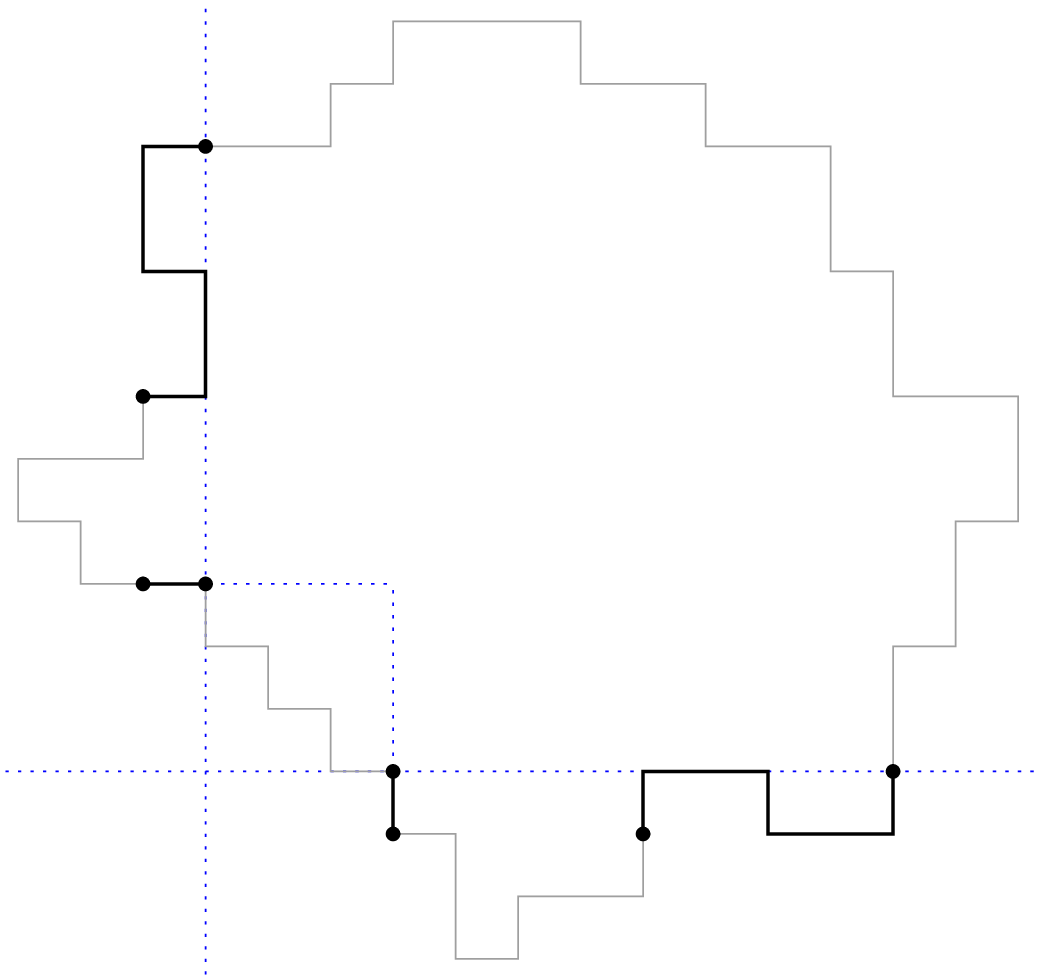}}
	\end{center}
	\caption{The form of 2-convex polygons with indents in different directions
	on adjacent sides.}
	\label{Fig:case8}
\end{figure}

The final case is enumerated in a similar fashion to the previous case. The
sub-classes defined by the relative heights of the indents are shown in
figure~\ref{Fig:case8}. The indents need not lie along the MBR and
thus the factors in the same quadrant as an indent must be of height
or width at least two.

There is, however, one special case to consider, depicted in part (e) of the
figure, namely when the  polygon does not enter the third quadrant.
This is an interesting case because
the horizontal and vertical joins (over $s$ and $t$ respectively) must be done
{\em simultaneously}. We cannot proceed is the usual way
and are forced to re-evaluate our approach and come up with a new way to evaluate
such constructions.

We solve this problem by generating the factor in the first quadrant as a
staircase factor and wrapping it along the factorisation lines to form a convex
factor. This is achieved by adding double-bonds joined to
the staircase factor next to the indent factors, so that it
wraps back along the double bonds, forming a convex polygon.
The length of the horizontal (vertical) join is counted by $s$ ($t$). 
Each horizontal double-bond has two steps, contributing
$s^2$ to the weight of the join, as well as a weight of $1/x$ for the fact that
it wraps the staircase back, reducing its width by one.
This means that the expression for the bottom factor in the join includes the
term $s^2/(x-s^2)$.

Other expressions including $s^2$ rather than just
$s$ have been evaluated preciously by simplifying the other side of the join
and re-expressing it in terms of derivatives. Now what is on the other
side of the join is the staircase factor, which cannot be simplified as it
has to be joined to both the bottom and left factors at the same time. We are
therefore forced to simplify the terms including $s^2$ or $t^2$. However since 
 $ s^2/(x-s^2) = E_x[s/(x-s)] $,  the
expression inside the $E$ operator can be changed
into a form we can evaluate simply.

The remaining cases are obtained using the above methodology by wrapping a
directed walk in the first quadrant. When the top indent is above the bottom one
the unimodal factor in the third quadrant must be of width and height
at least two .

\subsection{The 2-convex generating function}

In the above derivations, the direction and position of one of the indents was
chosen arbitrarily, such that the direction and position of the other determined
which sub-class the polygon belongs to. Therefore, when the two indents are in
different directions, the cardinality of the set of possible combinations of
direction and location for the fixed indent is four. When the indents are in the
same direction, both directions are enumerated by the generating functions.  If
the indents are on the same edge, the cardinality is two. If they are opposite,
it is one.  We therefore obtain the generating function for 2-convex
polygons by summing the results, multiplying each term by the cardinality of the
class it counts.  This gives the following generating function for 2-convex
polygons (available in a Maple/Mathematica friendly format
at  http://www.ms.unimelb.edu.au/\~{}iwan/polygons/series/2-convex-anisotropic.txt):
\begin{eqnarray} \nonumber
	\frac{-4A^{(c)}_2}{(1-x)^3 x^2 (1-y)^3 y^2 \Delta^{7/2}} \\
	- \frac{B^{(c)}_2}{(1-x)^7 x^2 (1-y)^7 y^2 ((1-x)^2-y)^3 ((1-y)^2-x)^3 (1-x-y) \Delta^4},
\end{eqnarray}
{\small
$
\mbox{\normalsize{where }} A^{(c)}_2 = 
\ (1-x)^{11} x^4
- 3 (1-x)^9 x^4 (5 - 2 x + x^2) y
+ (1-x)^7 x^2 (4 - 12 x + 103 x^2 - 79 x^3 + 31 x^4 - 11 x^5 + 3 x^6) y^2
- (1-x)^5 x^2 (40 - 124 x + 455 x^2 - 533 x^3 + 264 x^4 - 69 x^5 - 21 x^6 - 3 x^7 + x^8) y^3
+ (1-x)^3 (1 - 12 x + 232 x^2 - 742 x^3 + 1696 x^4 - 2297 x^5 + 1626 x^6 - 539 x^7 + 133 x^8 
+ 3 x^9 - 56 x^{10} + 3 x^{11}) y^4 + (1-x)^2 (-11 + 119 x - 943 x^2 + 2443 x^3 - 4014 x^4 
+ 4513 x^5 - 3054 x^6 + 867 x^7 - 58 x^8 + 221 x^9 - 137 x^{10} - 11 x^{11} + x^{12}) y^5
+ (1-x) (55 - 542 x + 2765 x^2 - 6154 x^3 + 8193 x^4 - 7901 x^5 + 5521 x^6 - 2140 x^7 - 284 x^8 
+ 430 x^9 + 19 x^{10} - 81 x^{11} + 7 x^{12}) y^6 + 
(-165 + 1503 x - 5996 x^2 + 11929 x^3 - 14004 x^4 + 11488 x^5 -
7661 x^6 + 4474 x^7 - 1456 x^8 - 506 x^9 + 504 x^{10} - 160 x^{11} + 18 x^{12}) y^7
+ (330 - 2502 x + 7381 x^2 - 10693 x^3 + 8925 x^4 - 4846 x^5 + 1856 x^6 - 1456 x^7 
+ 1364 x^8 - 328 x^9 + 103 x^{10} - 22 x^{11}) y^8
+ (-462 + 2898 x - 6353 x^2 + 6348 x^3 - 3639 x^4 + 1204 x^5 + 714 x^6 - 506 x^7 
- 328 x^8 - 32 x^9 + 10 x^{10}) y^9
+ (462 - 2394 x + 3925 x^2 - 2349 x^3 + 873 x^4 - 637 x^5 - 411 x^6 + 504 x^7 + 103 x^8 + 10 x^9) y^{10}
+ (-330 + 1422 x - 1797 x^2 + 434 x^3 + 47 x^4 + 484 x^5 - 100 x^6 - 160 x^7 - 22 x^8) y^{11}
+ (165 - 603 x + 632 x^2 - 4 x^3 - 180 x^4 - 114 x^5 + 88 x^6 + 18 x^7) y^{12}
- (1-x)^2 (55 - 67 x - 18 x^2 + 27 x^3 + 7 x^4) y^{13}
- (1-x)^3 (-11 + x^2) y^{14}
- (1-x)^3 y^{15}
$
\\[0.5\baselineskip]
\noindent
$
\mbox{\normalsize{and }} B^{(c)}_2 = 
\ 4 (1-x)^{26} x^4
- 4 (1-x)^{24} x^4 (30 - 20 x + 7 x^2) y
+ 4 (1-x)^{22} x^2 (4 - 12 x + 433 x^2 - 580 x^3 + 368 x^4 - 128 x^5 + 21 x^6) y^2
- 4 (1-x)^{20} x^2 (100 - 360 x + 4214 x^2 - 8192 x^3 + 7616 x^4 - 4297 x^5 
+ 1528 x^6 - 350 x^7 + 35 x^8) y^3
+ 4 (1-x)^{18} (1 - 12 x + 1252 x^2 - 5226 x^3 + 32426 x^4 - 76437 x^5 + 93156 x^6 
- 70591 x^7 + 36061 x^8 - 12917 x^9 + 3217 x^{10} - 532 x^{11} + 35 x^{12}) y^4
- (1-x)^{16} (104 - 1296 x + 42936 x^2 - 198752 x^3 + 859020 x^4 - 2164230 x^5 + 3180634 x^6 
- 3007804 x^7 + 1964426 x^8 - 927123 x^9 + 322834 x^{10} - 83390 x^{11} + 15842 x^{12} 
- 1961 x^{13} + 84 x^{14}) y^5
+ (1-x)^{14} (1300 - 16868 x + 289112 x^2 - 1412308 x^3 + 5065708 x^4 - 12664938 x^5 
+ 20907814 x^6 - 23411754 x^7 + 18540432 x^8 - 10767204 x^9 + 4716295 x^{10} - 1577652 x^{11} 
+ 403402 x^{12} - 79358 x^{13} + 12073 x^{14} - 1114 x^{15} + 28 x^{16}) y^6
- (1-x)^{12} (10400 - 140832 x + 1632176 x^2 - 8077604 x^3 + 26714752 x^4 - 64771716 x^5 
+ 113604188 x^6 - 143655868 x^7 + 132836186 x^8 - 91595522 x^9 + 48102042 x^{10} - 19641626 x^{11} 
+ 6301753 x^{12} - 1574826 x^{13} + 304364 x^{14} - 47486 x^{15} +6143 x^{16} - 400 x^{17} 
+ 4 x^{18}) y^7
- (1-x)^{10} (-59800 + 846768 x -7892920 x^2 + 38780872 x^3 - 125239164 x^4 + 297723354 x^5 
- 537008734 x^6 + 736108516 x^7 - 767512896 x^8 + 611776002 x^9 - 375779397 x^{10} 
+ 180183183 x^{11} - 68757620 x^{12} + 21225217 x^{13} - 5236463 x^{14} + 983039 x^{15} 
- 136742 x^{16} + 16579 x^{17} - 1952 x^{18} + 78 x^{19}) y^8
+ (1-x)^8 (-263120 + 3902272 x - 32692880 x^2 + 159495728 x^3 - 518222904 x^4 + 1235851899 x^5 
- 2273834564 x^6 + 3286427610 x^7 - 3744105010 x^8 + 3360169186 x^9 - 2372247082 x^{10} 
+ 1316981454 x^{11} - 578123374 x^{12} + 204616292 x^{13} - 60356681 x^{14} +  15025753 x^{15} 
- 2968864 x^{16} + 408498 x^{17} - 36191 x^{18} + 3653 x^{19} - 450 x^{20} + 7 x^{21}) y^9
+ (1-x)^6 (920920 - 14323848 x + 115445704 x^2 - 565821908 x^3 + 1881114592 x^4 - 4595055622 x^5 
+ 8689038821 x^6 - 13096803662 x^7 + 15935786206 x^8 - 15697620474 x^9 + 12485332431 x^{10} 
- 7965917640 x^{11} + 4040636749 x^{12} - 1619517767 x^{13} + 518745262 x^{14} - 140030433 x^{15} 
+ 34601374 x^{16} - 7726101 x^{17} + 1290999 x^{18} - 115181 x^{19} + 1171 x^{20} - 68 x^{21} 
+ 43 x^{22}) y^{10}
+ (1-x)^4 (-2631200 + 42965472 x - 346781248 x^2 + 1732729172 x^3 - 5960460908 x^4 + 15133786878 x^5 
- 29768776966 x^6 + 46871962702 x^7 - 60257832295 x^8 + 63871376160 x^9 - 55950160750 x^{10} 
+ 40346497476 x^{11} - 23700589523 x^{12} + 11139364197 x^{13} -  4087427426 x^{14} + 1145559467 x^{15} 
- 251575315 x^{16} + 52651169 x^{17} - 13478667 x^{18} + 3310337 x^{19} - 518486 x^{20} + 32936 x^{21} 
+ 909 x^{22} + 2 x^{23} + 3 x^{24}) y^{11}
- (1-x)^3 (-6249100 + 100966228 x - 786008036 x^2 + 3790621220 x^3 - 12634194004 x^4 
+ 31120468774 x^5 - 59333283779 x^6 + 90481426644 x^7 - 112836337344 x^8 + 116619494723 x^9 
- 100526269572 x^{10} + 72257250633 x^{11} -  42935196678 x^{12} + 20646926051 x^{13} 
- 7717297655 x^{14} + 2081066620 x^{15} - 344551302 x^{16} + 19708982 x^{17} - 190245 x^{18} 
+ 2604586 x^{19} - 1089925 x^{20} + 167990 x^{21} - 5325 x^{22} - 819 x^{23} + 5 x^{24}) y^{12}
- (1-x)^2 (12498200 - 200432672 x + 1523814896 x^2 - 7173568104 x^3 + 23394202064 x^4 
- 56432397385 x^5 + 105258925336 x^6 - 156756419659 x^7 + 190737577732 x^8 - 192647570536 x^9 
+ 163201860685 x^{10} - 116641434399 x^{11} + 70265808073 x^{12} - 35189781755 x^{13} 
+ 14129858643 x^{14} - 4206300270 x^{15} + 752961328 x^{16} + 2762280 x^{17} -41077210 x^{18} 
+ 8781488 x^{19} + 172960 x^{20} - 372435 x^{21} + 59206 x^{22} - 2116 x^{23} - 89 x^{24} 
+ 11 x^{25}) y^{13}
- (1-x) (-21246940 + 339297288 x - 2543202516 x^2 + 11796257160 x^3 - 37956243060 x^4 
+ 90389633588 x^5 - 166278246916 x^6 + 243693863948 x^7 - 291099521233 x^8 + 288163147560 x^9 
- 239402103854 x^{10} + 168782468152 x^{11} - 101938448323 x^{12} + 52836781903 x^{13} 
- 23089813175 x^{14} + 8076150924 x^{15} - 2006509792 x^{16} + 229306444 x^{17} 
+ 52941397 x^{18} - 29443771 x^{19} + 5233574 x^{20} - 137308 x^{21} - 96084 x^{22} 
+ 13023 x^{23} - 102 x^{24} + 33 x^{25}) y^{14}
+ (-30904640 + 493047872 x - 3671274848 x^2 + 16909765768 x^3 - 54076869168 x^4 
+ 128034624944 x^5 - 233957204524 x^6 + 339860131708 x^7 - 401161475162 x^8 
+ 390960587591 x^9 - 318456383150 x^{10} + 219480309921 x^{11} - 130108934416 x^{12} 
+ 67655799311 x^{13} - 31165964099 x^{14} + 12455845878 x^{15} - 4042728520 x^{16} 
+ 929084638 x^{17} - 93195098 x^{18} - 24059794 x^{19} +  11705293 x^{20} - 1983652 x^{21} 
+ 86582 x^{22} + 23547 x^{23} - 3151 x^{24} - 68 x^{25} + 5 x^{26}) y^{15}
+ (38630800 - 579105408 x + 4018988224 x^2 - 17153382704 x^3 + 50579891852 x^4 
- 109880934764 x^5 + 183275392200 x^6 - 241692656260 x^7 + 257400864147 x^8 
- 224539648706 x^9 + 161936707748 x^{10} - 97616424050 x^{11} + 50386570178 x^{12} 
- 23295420511 x^{13} + 10082660716 x^{14} - 4042728520 x^{15} + 1381926138 x^{16} 
- 352763968 x^{17} + 51556603 x^{18} + 1631397 x^{19} - 2559687 x^{20} + 547928 x^{21} 
- 54439 x^{22} - 611 x^{23} + 432 x^{24} - 31 x^{25}) y^{16}
+ (-41602400 + 585761792 x - 3791089760 x^2 + 15002579968 x^3 - 40793586584 x^4 
+ 81262920502 x^5 - 123537868000 x^6 + 147559428776 x^7 - 141331365760 x^8 
+ 109736530568 x^9 - 69193200719 x^{10} + 35421383132 x^{11} - 15013860403 x^{12} 
+  5709460646 x^{13} - 2235816236 x^{14} + 929084638 x^{15} - 352763968 x^{16} 
+ 103409188 x^{17} - 19757856 x^{18} + 1570918 x^{19} + 189358 x^{20} - 71908 x^{21} 
+ 11883 x^{22} - 268 x^{23} + 77 x^{24}) y^{17}
+ (38630800 - 510758608 x + 3084598320 x^2 - 11324232712 x^3 + 28395037220 x^4 
- 51826962892 x^5 + 71682098798 x^6 - 77353955608 x^7 + 66429125667 x^8 
- 45706907743 x^9 + 24906199156 x^{10} - 10403200527 x^{11} + 3174037717 x^{12} 
- 706359558 x^{13} + 176365047 x^{14} - 93195098 x^{15} + 51556603 x^{16} 
- 19757856 x^{17} + 4723092 x^{18} - 605060 x^{19} + 49569 x^{20} - 3395 x^{21} 
- 1253 x^{22} - 143 x^{23}) y^{18} 
+ (-30904640 + 383692992 x - 2163792656 x^2 + 7374341864 x^3 - 17053687520 x^4 
+ 28497776084 x^5 - 35790327136 x^6 + 34800116028 x^7 - 26729224736 x^8 
+ 16280988977 x^9 - 7648226004 x^{10} + 2524992746 x^{11} - 406853569 x^{12} 
- 93698188 x^{13} + 82385168 x^{14} - 24059794 x^{15} + 1631397 x^{16} + 1570918 x^{17} 
- 605060 x^{18} + 66554 x^{19} - 4454 x^{20} + 1819 x^{21} + 202 x^{22}) y^{19}
+ (21246940 - 247772008 x + 1305971612 x^2 - 4135556872 x^3 + 8823936024 x^4 
- 13492690928 x^5 + 15356870874 x^6 - 13416475166 x^7 +  9196907145 x^8 
- 4978222336 x^9 + 2052527732 x^{10} - 556811827 x^{11} + 29183400 x^{12} 
+ 58467226 x^{13} - 34677345 x^{14} + 11705293 x^{15} - 2559687 x^{16} + 189358 x^{17} 
+ 49569 x^{18} - 4454 x^{19} - 1368 x^{20} - 110 x^{21}) y^{20}
+ (-12498200 + 136980272 x - 675634696 x^2 + 1990453776 x^3 - 3921712444 x^4 
+ 5486790816 x^5 - 5649513210 x^6 + 4421089360 x^7 - 2695886222 x^8 + 1305174655 x^9 
- 491084981 x^{10} + 128534739 x^{11} - 11441768 x^{12} - 8063133 x^{13} 
+ 5370882 x^{14} - 1983652 x^{15} + 547928 x^{16} - 71908 x^{17} - 3395 x^{18} 
+ 1819 x^{19} - 110 x^{20}) y^{21}
+ (6249100 - 64433028 x + 297851620 x^2 - 817733112 x^3 + 1489860872 x^4 
- 1908408554 x^5 + 1775572130 x^6 - 1240442976 x^7 + 667928313 x^8 - 289317379 x^9 
+ 102644601 x^{10} - 29961766 x^{11} + 6383656 x^{12} - 977036 x^{13} - 41224 x^{14} 
+ 86582 x^{15} - 54439 x^{16} + 11883 x^{17} - 1253 x^{18} + 202 x^{19}) y^{22}
+ (-2631200 + 25559072 x - 110940720 x^2 + 284446548 x^3 - 480234432 x^4 
+ 564112724 x^5 - 474248936 x^6 + 294904852 x^7 - 138225516 x^8 + 52522411 x^9 
- 17224508 x^{10} + 5578263 x^{11} - 1609051 x^{12} + 492963 x^{13} - 109107 x^{14} 
+ 23547 x^{15} - 611 x^{16} - 268 x^{17} - 143 x^{18}) y^{23}
+ (920920 - 8444128 x + 34492144 x^2 - 82799336 x^3 + 129893340 x^4 - 140307634 x^5 
+ 106798702 x^6 - 59095890 x^7 + 23639299 x^8 - 7486656 x^9 + 2001655 x^{10} 
- 644781 x^{11} + 181503 x^{12} - 63349 x^{13} + 13125 x^{14} - 3151 x^{15} 
+ 432 x^{16} + 77 x^{17}) y^{24}
+ (-263120 + 2283072 x - 8798288 x^2 + 19829008 x^3 - 28995448 x^4 + 28906719 x^5 
- 20012600 x^6 + 9913898 x^7 - 3347134 x^8 +  826000 x^9 - 124087 x^{10} + 29300 x^{11} 
- 2853 x^{12} + 1927 x^{13} - 135 x^{14} - 68 x^{15} - 31 x^{16}) y^{25}
+ (59800 - 492568 x + 1795744 x^2 - 3809804 x^3 + 5207920 x^4 - 4809154 x^5 
+ 3046581 x^6 - 1369614 x^7 + 399732 x^8 - 78407 x^9 +  2224 x^{10} + 919 x^{11} 
- 834 x^{12} + 111 x^{13} + 33 x^{14} + 5 x^{15}) y^{26}
+ (-10400 + 81632 x - 282336 x^2 + 565028 x^3 - 723116 x^4 + 619222 x^5 - 359946 x^6 
+ 148482 x^7 - 39609 x^8 + 7449 x^9 - 326 x^{10} - 10 x^{11} + 5 x^{12} - 11 x^{13}) y^{27}
+ (1300 - 9768 x + 32140 x^2 - 60712 x^3 + 72592 x^4 - 57298 x^5 + 30217 x^6 - 11207 x^7 
+ 2732 x^8 - 506 x^9 + 43 x^{10} + 3 x^{11}) y^{28}
+ (-104 + 752 x - 2360 x^2 + 4200 x^3 - 4648 x^4 + 3305 x^5 - 1506 x^6 + 448 x^7 
- 78 x^8 + 7 x^9) y^{29}
+ 4 (1-x)^7 y^{30}.
$
}

\section{Summary and outlook}

Following Lin's approach to enumerating 1-convex polygons, one can factorise
almost-convex polygons by extending lines along the base of all indents. Then,
using a `divide and conquer' approach, it is possible to then enumerate the
various sub-classes.

By looking at the form of the various factors, it is possible to guess what the
form of the resulting generating functions will be. One can then obtain the
generating functions by directly enumerating the series to a sufficiently high
order and then solving the set of linear equations corresponding to the presumed
form.

We presented some techniques that are invaluable in enumerating the factors
exactly. The most important technique is that of `wrapping', which allows the
generation of quite complex objects out of simply enumerable components. It is easily
implemented when using an inclusion-exclusion approach to enumerating unimodal
factors, which made it essential in enumerating 2-convex polygons. These
techniques have allowed us to reduce the enumeration of 1-convex polygons, as
well as many sub-classes of 2-convex polygons, to a single, combinatorially
interpretable expression.

Going forward, it is not realistic to factorise almost-convex polygons
for high concavity indices as we have done here, as there will be an exponential 
growth in the number of cases to evaluate. It would be more sensible
to restrict the size of each indentation first, and then generalise these cases.
Eventually, defining operators that can add more and more complex indentations
in the side of convex polygons and looking at the effect on the asymptotic
growth of their number seems to be the most appropriate path to understand
 how convex polygons become general SAPs.

\ack

We gratefully acknowledge financial support from the Australian Research Council.

\section*{References}


\end{document}